\documentclass[12pt,a4paper]{article}

\usepackage{pdflscape}
\usepackage{amsmath}
\usepackage{amssymb}
\usepackage{latexsym}
\usepackage{srcltx}
\usepackage{graphics}
\usepackage{color}
\usepackage{epsfig}
\usepackage{color}
\usepackage[ruled,vlined]{algorithm2e}
\usepackage{GGGraphics}

\newcommand{\co}{{\mathbb C}}

\newcommand{\R}{{\mathbb R}}
\newcommand{\re}{{\mathbb R}}
\newcommand{\n}{{\mathbb N}}

\renewcommand{\P}{\Pi}
\newcommand{\z}{{\mathbb Z}}
\newcommand{\cA}{{\mathcal{A}}}
\newcommand{\cC}{{\mathcal{C}}}

\newcommand{\cV}{{\mathcal{V}}}

\newcommand{\cB}{{\mathcal{B}}}

\newcommand{\cH}{{\mathcal{H}}}
\newcommand{\cS}{{\mathcal{S}}}
\newcommand{\cR}{{\mathcal{R}}}

\newcommand{\cP}{{\mathcal{P}}}

\newcommand{\pa}{{\mathbf{\alpha}}}

\newcommand{\tA}{{\tilde A}}

\newtheorem{theorem}{Theorem}

\newtheorem{lemma}{Lemma}
\newtheorem{cor}{Corollary}
\newtheorem{remark}{Remark}
\newtheorem{ex}{Example}
\newtheorem{defi}{Definition}

\newtheorem{ass}{Assumption}

\date{}

\author{Nicola Guglielmi
\thanks{Dipartimento di 
Matematica and DEWS,
University of L'Aquila, Italy
{e-mail: \tt\small guglielm@univaq.it }}
 and
Vladimir Yu.~Protasov
\thanks{Dept. of Mechanics and Mathematics, Moscow State University,
Vorobyovy Gory, 119992, Moscow, Russia,  {e-mail: \tt\small
v-protassov@yandex.ru}}}

\title{Invariant  polytopes
of linear operators\\ with applications to regularity of wavelets \\ and of subdivisions
\thanks{
The first author is supported by INdAM GNCS (Gruppo Nazionale di Calcolo Scientifico); the second author is supported by RFBR grants nos. 13-01-00642 and 14-01-00332, and by the grant of Dynasty foundation.}}

\begin{document}
\maketitle

\begin{abstract}

We generalize  the recent invariant  polytope algorithm for computing  the joint spectral radius and extend it to a wider class of matrix sets. This, in particular, makes the algorithm applicable to sets of matrices that have finitely  many spectrum maximizing products. A criterion of convergence of the algorithm is proved.

As an application we solve two challenging computational open problems. First we find the regularity of  the Butterfly subdivision scheme for various parameters $\omega$. In the ``most regular'' case $\omega = \frac{1}{16}$, we prove that the limit function has  H\"older exponent  $2$ and its derivative is ``almost Lipschitz'' with logarithmic factor~$2$. Second we compute the H\"older exponent of Daubechies wavelets of high order.
\smallskip

\smallskip


\end{abstract}
\bigskip

\section{Introduction}\label{s1}

The {\em joint spectral radius} of a set of matrices (or linear operators) originated  in early sixties with Rota and Strang~\cite{RS60} and found countless applications in functional analysis, dynamical systems, wavelets, combinatorics, number theory, automata, formal languages, etc. (see bibliography in~\cite{DL92, G, GP13, J09}). We focus on finite sets of matrices, although all the results are extended to arbitrary compact sets. If the converse is not stated, we assume a fixed basis in $\re^d$
and identify an operator with the corresponding matrix. Everywhere below $\cA = \{A_1, \ldots , A_m\}$ is an arbitrary family of $d\times d$-matrices, $\cA^k$ is the set of all $m^k$ products of $k$ matrices from~$\cA$ (with repetitions permitted). A product 
$\Pi \in \cA^k, \, n \in \n$,  is called {\em simple} if it is not a power of a shorter product.
\begin{defi}\label{d3}
The joint spectral radius of a family $\cA$ is {\rm \cite{RS60}}
\begin{equation}\label{jsr}
\rho(\cA)\ = \ \lim_{k\to \infty} \max_{A \in \cA^k} \|A\|^{1/k}\, .
\end{equation}
\end{defi}
This limit exists and does not depend on the matrix norm.
In case $m=1$ (i.e., $\cA = \{A_1\}$), according to Gelfand's theorem, the joint spectral radius is $\lim_{k\to \infty}\|A_1^k\|^{1/k}$, i.e., coincides with the usual spectral radius $\rho(A_1)$, which is the maximal modulus of eigenvalues of~$A_1$ (see, for instance~\cite{BW92,E95}).

The joint spectral radius has the following geometrical meaning: $\rho(\cA)$ is the infimum of numbers $\nu$ for which there is a norm $\|\cdot \|_{\nu}$ in $\re^d$ such that in the induced operator norm, we have $\|A_i\|_{\nu} \le \nu, \, i = 1, \ldots , m$.
In particular, $\rho(\cA) < 1$ if and only if all operators from~$\cA$ are contractions in some 
(common) norm. Thus, the joint spectral radius is the indicator of simultaneous contractivity of operators $A_1, \ldots , A_m$.

Another interpretation is due to the discrete dynamical system:
$$
x(k+1) = A(k) x(k)\, , \quad k = 0, 1, 2,  \ldots ,
$$
where each matrix $A(k)$ is chosen from $\cA$ independently for every $k$ and $x(0) \in \re^d \setminus \{0\}$.
Then $\rho(\cA)$ is the exponent of fastest possible growth of trajectories of the system:
the maximal upper limit
$\limsup\limits_{k \to \infty} \frac{\log \|x_k\|}{k}$ among all the trajectories $\{x_k\}_{k \ge 0}$
is equal to $\log \rho(\cA)$.
In particular, $\rho(\cA) < 1$ if and only if the system is stable, i.e., $x(k) \to 0$ as $k \to \infty$, for every
trajectory.

The problem of computing or estimating the joint spectral radius is notoriously hard.  This is natural in view of the negative complexity results of Blondel and Tsitsiklis~\cite{BT1, BT2}. Several methods of approximate computation  were elaborated in the literature~\cite{AS, BN, BNT, CH, DL92, G, GWZ05, PJ, P96, P97, PJB}.  They work well in low dimensions (mostly, not exceeding $5-8$). When the dimension growth, then ether
the estimation becomes rough or the the running time grows dramatically. In recent work~\cite{GP13} we derived the {\em invariant polytope algorithm} that finds the exact value of $\rho(\cA)$ for the vast majority of matrix families in dimensions up to $20$.  For sets of nonnegative matrices it works faster and finds the exact value in dimensions up to $100$ and higher. We conjectured in~\cite{GP13} that the set of matrix families~$\cA$ for which this algorithm finds $\rho(\cA)$ within finite time is of full Lebesgue measure. Several open problems from applications have been solved by using this method. However, it cannot handle one important case, which often emerges in practice:
the case of several {\em spectrum maximizing products}. In this paper we generalize this method making it applicable for a wider class of matrix families, including that special case. To formulate the problem we need some more facts and notation.

The following double inequality for the joint spectral radius $\rho = \rho(\cA)$ is well known:
\begin{equation}\label{ineq}
\max_{A \in \cA^k} [\rho(A)]^{1/k} \ \le \ \rho \ \le \    \max_{A \in \cA^\ell} \|A\|^{1/\ell}\, , \quad k,\ell \in \n\, .
\end{equation}
Moreover, both parts of this inequality converge to $\rho$ as $k \to \infty$. In fact,\linebreak
$\lim\limits_{k\to \infty} \max\limits_{A \in \cA^k} \|A\|^{1/k} = \rho$ (by definition) and
 $\limsup\limits_{k\to \infty} \max\limits_{A \in \cA^k} [\rho(A)]^{1/k} = \rho$ (see~\cite{BW92}). All algorithms of approximate computation
 of the joint spectral radius are based on this inequality. First, one finds the greatest value of 
$\max_{A \in \cA^k} [\rho(A)]^{1/k}$ (the left hand side of~(\ref{ineq})) among all reasonably small~$k$, then one minimizes
the value  $ \max_{A \in \cA^k} \|A\|^{1/k}$ (the right hand side of~(\ref{ineq}))  by choosing an appropriate matrix norm $\|\cdot \|$.
Thus maximizing the lower bound and minimizing the upper one we approximate the joint spectral radius.
Sometimes those two bounds meet each other, which gives the exact value of $\rho$. This happens when one finds the {\em spectrum maximizing product}~$A = \Pi \in \cA^k$ and the norm $\| \cdot \|$ for which both inequalities in~(\ref{ineq}) become equalities.
\begin{defi}\label{d5}
A simple product $\Pi \in \cA^n$ is called the spectrum maximizing product (s.m.p.) if the value $[\rho(\Pi)]^{1/n}$
is maximal among all products of matrices from~$\cA$ of all lengths~$n \in \n$.
\end{defi}
Let us remark that an s.m.p. maximizes the value $[\rho(\Pi)]^{1/n}$ among all products of our matrices,
not just among products of length~$n$.
Observe that for any s.m.p., $\Pi \in \cA^n$, we have $[\rho(\Pi)]^{1/n} = \rho(\cA)$. Indeed, from~(\ref{ineq}) it follows that
$[\rho(\Pi)]^{1/n} \le \rho(\cA)$. If this inequality is strict, then there are $k\in \n$
such that  $\max_{A \in \cA^k} [\rho(A)]^{1/k} > [\rho(\Pi)]^{1/n}$ (because of the convergence property),
which contradicts to the maximality of~$\Pi$. Thus, {\em to find the joint spectral radius it suffices to prove that a given product $\Pi \in \cA^n$ is an s.m.p.}.
The invariant polytope algorithm~\cite{GP13}  proves the s.m.p. property of a chosen product~$\Pi$ by recursive construction of a polytope 
(or more general $P \subset \co^d$, although here we consider for simplicity real polytopes) $P \subset \re^d$ such that $A_i P  \subset [\rho(\Pi)]^{1/n} P, \, i = 1, \ldots , m$.

In the Minkowski norm $\|\cdot \|_P$ generated in $\re^d$ by this polytope, we have
$\max_{A \in \cA} \|A\|_P \le [\rho(\Pi)]^{1/n}$ and $\|\cdot \|_P$ is said to be an
extremal norm for $\cA$. Applying~(\ref{ineq}) for $k  = n$ (left hand side inequality) and for $\ell = 1$ (right hand side) we conclude that $[\rho(\Pi)]^{1/n} = \rho(\cA)$.

Note that if a product $\Pi \in \cA^n$ is s.m.p., then so is each of $n$~its cyclic permutations.
If there are no other s.m.p., then we say that the s.m.p. is {\em unique} meaning that it is unique  up to cyclic permutations.

The disadvantage of the  polytope algorithm is that
it is guaranteed to work only if the family~$\cA$ has a unique s.m.p. Otherwise the algorithm may not be able to construct the desired polytope within finite time, even
if this exists.
The uniqueness of s.m.p. condition, although believed to be generic, is not satisfied in many practical cases.
For example, it happens often that several matrices  $A_i \in \cA$ are s.m.p.
(of length~$1$). The extension of our algorithm presented in this paper works with an arbitrary 
(finite) number of s.m.p.'s.
We prove the theoretical criterion of  convergence of the algorithm and apply it to solve
two long standing open problems: computing  the H\"older regularity of the Butterfly subdivision scheme (Section~\ref{s5}) and computing the regularity of Daubechies wavelets of high order (Section~\ref{s6}).

\section{Statement of the problem}\label{s2}

We begin with a short description of the invariant polytope algorithm from~\cite{GP13}
for computing the joint spectral radius and finding an extremal polytope norm of a given family $\cA = \{A_1, \ldots , A_m\}$. We make a usual assumption that the family is irreducible, i.e., its matrices do not have a common nontrivial invariant subspace. Recall that for a reducible family the  computation of the joint spectral radius is obtained by solving several problems of smaller dimensions~\cite{DL92}.
For a given set $M \subset \re^d$ we denote by ${\rm co}(M)$ the convex hull of $M$ and 
by ${\rm absco}(M) = {\rm co}\{M, -M\}$ the symmetrized convex hull. The sign $\asymp$ denotes as usual the asymptopic equivalence of two values (i.e., equivalence up to multiplication by a constant).
\smallskip

\noindent \textbf{The invariant polytope algorithm} (see \cite{P97,GWZ05,GP13}).
\smallskip

{\tt Initialization.} First, we fix some number ${\overline n}$ and find a simple
product~$\Pi = A_{d_n}\ldots A_{d_1}$ with the maximal value $[\rho(\Pi)]^{1/n}$ among all products of lengths $n \le {\overline n}$. We call this product a {\em candidate s.m.p.} and try to prove that it is actually an s.m.p. Denote $\rho_c = [\rho(\Pi)]^{1/n}$ and normalize all the matrices $A_i$ as $\tilde A_i = \rho_c^{-1}A_i$.
Thus we obtain the family~$\tilde A$ and the product $\tilde \Pi = \tilde A_{d_n}\ldots \tilde A_{d_1}$ such that $\rho(\tilde \Pi) = 1$. For the sake of simplicity we assume that the largest by modulo eigenvalue of~$\tilde \Pi$ is real, in which case it is $\pm 1$. We assume it is $1$, the case of $-1$ is considered in the same way. 
The eigenvector $v^{(1)}$ corresponding to this eigenvalue is called {\em leading eigenvector}.
The vectors $v^{(j)} = \tilde A_{d_{j-1}}\cdots \tilde A_{d_1}v^{(1)}, \, 
j= 2, \ldots , n$, are leading eigenvectors of cyclic permutations of~$\tilde \Pi$.
The set $\cH = \{v^{(1)}, \ldots , v^{(n)}\}$ is called {\em root}.
Then we construct a sequence of finite sets $\cV_i \subset \re^d$ and their
subsets~$\cR_i \subset \cV_i$ as follows:

{\tt Zero iteration}. We set $\cV_0 = \cR_0 = \cH$.

{\tt $k$-th iteration, $k \ge 1$}. We have finite set $\cV_{k-1}$ and its subset
$\cR_{k-1}$.
We set $\cV_{k} = \cV_{k-1},  \, \cR_k = \emptyset$ and for every $v \in \cR_{k-1}, \, \tilde A\in \tilde \cA$, check whether $\tilde Av$ is an {\em interior point} of
${\rm absco} (\cV_{k})$ (this is an LP problem).
If so, we omit this point and take the next pair $(v , \tilde \cA) \in \cR_{k-1}\times \tilde \cA$, otherwise we add $\tilde A v$ to $\cV_{k}$ and to $\cR_k$. When all pairs $(v, \tilde A)$ are exhausted, both $\cV_{k}$ and $\cR_k$ are constructed.
Let $P_k = {\rm absco}(\cV_k)$. We have
$$
\cV_{k}
= \, \cV_{k-1} \cup \cR_k\, , \quad
P_{k} = {\rm co}\,  \{\tilde A_1P_{k-1}, \ldots
  , \tilde A_m P_{k-1}\}\, .
$$

{\tt Termination}. The algorithm halts when $\cV_{k} = \cV_{k-1}$, i.e., $\cR_k = \emptyset$ (no new vertices are
added in the $k$-th iteration). In this case
$P_{k-1} = P_k$, and hence $P_{k-1}$ is an invariant polytope, $\Pi$ is an s.m.p., and $\, \rho(\cA)  = [\rho(\Pi)]^{1/n}$.
\hfill \textbf{End of the algorithm}.
\smallskip

Actually, the algorithm  works with the sets $\cV_k$  only, the polytopes~$P_k$ are needed to illustrate the idea. Thus, in each iteration of the algorithm, we construct a polytope $P_k \subset \re^d$, store all its vertices
in the set~$\cV_k$ and spot the set $\cR_k \subset \cV_k$ of newly appeared (after the previous iteration) vertices.
Every time we check whether $\tilde \cA P_k \subset P_k$. If so, then $P_k$ is an invariant polytope, 
$\|\tilde A_i\|_{P_k} \le 1$ for all $i$, where $\|\cdot \|_{P_k}$ is the Miknowski norm associated 
to the polytope~$P_k$, and $\Pi$ is an s.m.p. Otherwise, we update the sets $\cV_k$ and $\cR_k$ and continue.

If the algorithm terminates within finite time, then it proves that the
chosen candidate is indeed an s.m.p. and gives the corresponding polytope norm. Although  there are simple examples of matrix families
for which the algorithm does not terminate, we believe that such cases are rare in practice.
In fact, in all numerical experiments made with randomly generated matrices and with matrices from applications, the algorithm did terminate in finite time providing an invariant polytope. The only special case when it does not work is when there are several different s.m.p. (up to cyclic permutations). In this case the algorithm never converges as it follows from the criterion proved in~\cite{GP13}. The criterion uses the notion of dominant product which is a strengthening of the s.m.p. property. A product
$\Pi \in \cA^n$ is called {\em dominant} for the family $\cA$ if there is a constant
$\gamma < 1$ such that the spectral radius of each product of matrices
from the normalized family $\tilde \cA = [\rho(\Pi)]^{-1/n}\cA$, which is neither a power 
of~$\tilde \Pi$ nor one of its cyclic permutations, is smaller than~$\gamma$.
A dominant product is an s.m.p., but, in general, not vice versa.

\smallskip

\noindent \textbf{Theorem A}~\cite{GP13}. {\em For a given set of matrices and for a given initial product~$\Pi$,
the invariant polytope algorithm terminates within finite time if and only if~$\Pi$ is dominant and its
leading eigenvalue is unique and simple.}
\smallskip

Note that if there is another s.m.p., which is neither a power of~$\Pi$ nor of its cyclic permutation,
then $\Pi$ is not dominant. Therefore, from  Theorem~A we conclude
\begin{cor}\label{c10}
If a family~$\cA$ has more than one s.m.p., apart from taking powers or cyclic permutations, then, for every initial product, the invariant polytope algorithm  does not terminate within finite time.
\end{cor}
The problem occurs in the situation when a family has several s.m.p., although not generic, but possible in some relevant applications. Mostly those are s.m.p. of length
$1$, i.e., some of matrices of the family $\cA$ have the same spectral radius and dominate the others. This happens, for instance, for transition matrices of refinement equations and wavelets (see Sections~\ref{s5}, \ref{s6}).
In the next section we show that the algorithm can be modified and extended to families with finitely many spectrum maximizing products.

Let a family $\cA$ have $r$ candidate s.m.p.'s $\Pi_1, \ldots , \Pi_r, \, r \ge 2$.
These products are assumed to be simple and different (up to cyclic permutations). Denote by $n_i$ the length of $\Pi_i$. Thus, $[\rho(\Pi_1)]^{1/n_1} = \ldots = [\rho(\Pi_r)]^{1/n_r} = \rho_c$.
Let $\tilde \cA = \rho_c^{-1} \cA$ be the normalized family,
$v_i$ be a leading eigenvector
of $\tilde \Pi_i$ (it is assumed to be real) and $\cH_i = \{v_i^{(1)}, \ldots , v_i^{(n_i)}\}$
be the corresponding roots, $i=1 \ldots , r$. The first idea is to start constructing
 the invariant polytope with all roots simultaneously, i.e., with the initial set of vertices
 $$
 \cV_0 \ =\ \cup_{i=1}^{\, r} \cH_i \ = \  \Bigl\{\ v_i^{(k)}  \ \Bigl| \ k = 1, \ldots , \ n_i\, , \ i = 1, \ldots , r\,
 \Bigr\}\, .
 $$
 However, this algorithm may fail to converge as the following example demonstrates.
 \begin{ex}\label{ex10}{\em Let $A_1$ and $A_2$ be operators in~$\re^2$:
 $A_1$ is a contraction with factor~$\frac12$ towards the $OX$ axis along the vector
$(1, 4)^T$, $A_2$ is a contraction with factor~$\frac12$ towards the $OY$  axis along the vector $(1, -2)^T$.
The matrices of these operators are
$$
A_1 \ = \
\left(
\begin{array}{rr}
1 & -\frac18 \\
0 & \frac12
\end{array}
\right) \
; \qquad
A_1 \ = \
\left(
\begin{array}{rr}
\frac12 & 0 \\
1 & 1
\end{array}
\right)
$$
Clearly, both $A_1$ and $A_2$ have a unique simple leading eigenvalue~$1$; $v_1 = (1,0)^T$ is the leading eigenvector of~$A_1$ and $v_2 = (0,1)^T$ is the leading eigenvector of~$A_2$.

The algorithm with two candidate s.m.p.'s
$\Pi_1 = A_1, \Pi_2 = A_2$ and with initial vertices
$\cV_{0} = \{v_1, v_2\}$ does not converge. Indeed, the set ${\rm absco} (\cV_k)$
has an extreme point~$A_2^kv_1$, which tends to the point $2v_2$ as $k \to \infty$,
but never reaches  it.

On the other hand, the same algorithm with initial vertices $\cV_0 = \{v_1, 3v_2\}$
terminates immediately after the first iteration with the invariant polytope
(rhombus)~$P = {\rm absco} (\cV_0) = {\rm co} \{\pm v_1, \pm 3v_2\}$.
Indeed, one can easily check that $A_iP \subset P, \, i = 1,2$. Therefore,
$A_1$ and $A_2$ are both s.m.p. and $\rho(\cA) = 1$.}
\end{ex}

This example shows that if the algorithm does not converge with the
leading eigenvectors $v_1, \ldots , v_r$ (or with the roots
$\cH_1, \ldots , \cH_r$), it, nevertheless, may converge if one multiplies these
eigenvectors (or the roots) by some numbers $\alpha_1, \ldots , \alpha_r$.
In Example~\ref{ex10} we have $\alpha_1 = 1, \alpha_2 = 3$. We call a vector of positive multipliers 
$\pa = (\alpha_1, \ldots , \alpha_r)$ {\em balancing vector}.

Thus, if a family has several candidate s.m.p.'s, then one can {\em balance} its
leading eigenvectors (or its roots) i.e., multiply them by the entries $\{\alpha_i\}_{i=1}^r$ of 
some balancing vector~$\pa > 0$  and start
the invariant polytope algorithm. In the next section we prove that the balancing
vector~$\pa$ for which the algorithm converges does exist and can be efficiently found, provided all the products $\Pi_1, \ldots , \Pi_r$ are dominant (meaning the natural extension of dominance from a single product to a set of products).
Thus, the corresponding extension of the invariant polytope algorithm is given by
Algorithm \ref{algoP}.

\begin{algorithm}
\KwData{$\cA = \{ A_1,\ldots,A_m \}, \ k_{\max}$}
\KwResult{The invariant polytope $P$, spectrum maximizing products, the joint spectral 
radius~$\rho(\cA)$}
\Begin{
\nl Compute a set of candidate spectrum maximizing products $\Pi_1,\ldots,\Pi_r$\;
\nl Set $\rho_c := \rho(\P_1)^{1/n_1}$ and $\tilde{\cA} := \rho_c^{-1} \cA$\;
\nl Compute $v_1,\ldots,v_r$, leading eigenvectors of $\tilde \Pi_1,\ldots,\tilde \Pi_r$ with $\| v_j \| = 1$ for all $j$\;
\nl Form the roots $\cH_1, \ldots , \cH_r$\;
\nl Provide the
 positive scaling factors $\alpha_1,\ldots,\alpha_r \le 1$\;
\nl Set $\cV_0 := \{\alpha_j \cH_j\}_{j=1}^{r}\, \ \cR_0 := \cV_0$\;
\nl Set $k=1$\;
\nl Set $E=0$\;
\nl \While{$E=0$ \ {\rm and} \ $k \le k_{\max}$}{
\nl Set $\cV_{k} = \cV_{k-1}, \, \cR_{k} = \emptyset$\;
\nl \For{$v \in \cR_{k-1}, \, \mbox{{\rm \textbf{and for}}} \ i = 1, \ldots , m$}{
\nl \eIf{$\tilde A_iv \in {\rm int} ({\rm absco}(\cV_{k}))$}{
\nl Leave $\cV_k, \cR_k$ as they are\;}{
\nl Set $\cV_{k} := \cV_{k}\cup \{\tilde A_i v\}, \, \cR_{k} := \cR_{k}\cup \{\tilde A_i v\}$\;}}
\nl \eIf{$\cR_{k} = \emptyset$}{
\nl Set $E=1$ (the algorithm halts) \; }{
\nl Set $k:= k+1$ \;}
}
\nl \eIf{$E=1$}{
\nl \Return{$P := {\rm absco}(\cV_k) \ $ is an invariant polytope;\\
$\hspace{14mm} \tilde \Pi_1, \ldots \tilde \Pi_r$ are s.m.p.;\\
$\hspace{14mm} \rho(\cA) = \rho_c$  is the joint spectral radius\;}}
{{\bf print} {\rm Maximum number of iterations reached}\;}
}
\caption{The invariant  polytope algorithm extension \label{algoP}}
\end{algorithm}

A crucial point in Algorithm \ref{algoP} is Step {\bf 5}, where suitable scaling factors 
$\alpha_1, \ldots , \alpha_r$ have
to be given. Then Algorithm~\ref{algoP} essentially repeats the
invariant polytope algorithm replacing the root $\cH$ by the union of roots 
$\alpha_1\cH_1, \ldots , \alpha_r\cH_r$.
Finding the scaling factors that provide the convergence of the algorithm is a nontrivial problem.
A proper methodology to compute them in an optimal way (in other words, to balance the leading eigenvectors) is derived  in the next section.
\begin{remark}\label{r40}
{\em Till now we have always assumed that the leading eigenvalue is real. This is not a restriction, because the invariant polytope algorithm is generalized to the complex case as well~(Algorithm~C in \cite{GP13,GWZ05,GZ07}). For the sake of simplicity, in this paper we consider only the case of real leading eigenvalue.}
\end{remark}

\section{Balancing leading eigenvectors. The main results}\label{s3}

\subsection{Definitions, notation, and auxiliary facts}\label{s3.1}

Let $\Pi_1, \ldots , \Pi_r$ be some products of matrices from~$\cA$ of lengths $n_1, \ldots , n_r$ respectively.
They are assumed to be simple and different up to cyclic permutations.
We also assume that all those products are candidates s.m.p.'s, in particular, 
$[\rho(\Pi_1)]^{1/n_1} = \ldots = [\rho(\Pi_r)]^{1/{n_r}} = \rho_c$.
We set $\tilde \cA = \rho_{c}^{-1} \cA$ and denote as $M$ the supremum of norms of all products of matrices from $\tilde \cA$.
Since $\tilde \cA$ is irreducible and $\rho (\tilde \cA) = 1$, this supremum is finite~\cite{BW92}.
By $\cA^* =  \{ A_1^*,\ldots ,  A_m^*\}$ we denote the family of adjoint matrices, the definition 
of~$\tilde \cA^*$ is analogous.

To each product $\tilde A_{b_n}\ldots \tilde A_{b_1}$ we associate the 
word $b_n\ldots b_1$ of the alphabet $\{1, \ldots , m\}$. By $ab$ we denote concatenation of words $a$ and $b$, in particular, $a^n = a\ldots a$ ($n$ times); the length of the word $a$ is denoted by $|a|$.
A word is {\em simple}  if it is not a power of a shorter word. In the sequel we identify
words with corresponding products of matrices from~$\tilde \cA$.
\begin{ass}
We now make the main assumption:  each product $\tilde \Pi_i$ has a real leading eigenvalue (Remark~\ref{r40}), which is either $1$ or $-1$ in this case.
For the sake of simplicity, we assume that all $\lambda_i = 1$ (the case $\lambda_i=-1$ is considered in the same way).
We denote by $v_i$ the corresponding leading eigenvector of~$\Pi_i$
(one of them, if it is not unique).
\end{ass}
Clearly, the corresponding adjoint matrix $\tilde \Pi^*_i$  also has a leading eigenvalue $1$, and a real leading eigenvector $v_i^*$. If the the leading eigenvalue is unique and simple, then  $\bigl(v_i^*, v_i \bigr) \ne 0$. In this case we normalize the adjoint leading eigenvector~$v_i^*$ by the condition  $\bigl(v_i^*, v_i \bigr) = 1$.

Take arbitrary $i = 1, \ldots , r$ and consider the product 
$\tilde \Pi_i = \tilde A_{d_n}\cdots \tilde A_{d_1}$, where $n = n_i$ and $d_s = d_s^{(i)}$
(for the sake of simplicity we omit the indices~$i$).  The vectors $v_i^{(1)}  = v_i, \, v_i^{(2)} = \tilde A_{d_1}v_i, \, \ldots \, , \, v_i^{(n)} = \tilde A_{d_{n-1}} \cdots \tilde A_{d_1}v_i$ are the leading eigenvectors of cyclic permutations of the product~$\Pi_i$.
The set $\cH_i = \{v_i^{(1)}, \ldots , v_i^{(n)}\}$ is a root of the tree from which the polytope algorithm starts.
Let $P_{i, k} = {\rm absco}\, \bigr\{\, \tilde \cA^p \, \cH_i \ \bigr| \ p = 0, \ldots , k\, \bigr\}$ be the polytope produced after the $k$-th iteration of the algorithm
started with the product $\Pi_i$, or, which is the same, with the root~$\cH_i$.
This polytope is a symmetrized convex hull of the set $\cV_{i, k}$ of all alive vertices of the tree on the first $k$ levels. In particular,
$\cV_{i, 0} = \cH_{i}$ and $P_{i, 0} = {\rm absco } (\cH_{i})$. We denote by $P_{i, \infty}$ and $\cV_{i, \infty}$ the union of the corresponding sets
for all $k \in \n \cup \{0\}$. If Algorithm~\ref{algoP} terminates within finite time, then $\cV_{i, \infty}$ is finite and $P_{i, \infty}$
is a polytope. If $\tilde \Pi_i$ is an s.m.p., i.e, if $\rho (\tilde \Pi_i) = 1$, then, by the irreducibility, the set $P_{i, \infty}$ is bounded~\cite{BW92}.
\smallskip

Now assume that all $\Pi_i$ have unique simple leading eigenvalues, in which case all the adjoint leading eigenvectors~$v_i^*$ are normalized by the condition $(v_i^*, v_i) = 1$. 
For an arbitrary pair $(i, j) \in \{1, \ldots , r\}^2$, and for arbitrary $k \in \{0\}\cup \n \cup \{\infty\}$, we denote
\begin{equation}\label{qij}
q_{ij}^{(k)} \ =
\ \sup_{z \in \cV_{i, k}}\, \bigl|\, \bigl(v^*_j\, , \, z\bigr)\, \bigr|,
 \end{equation}
Thus, $q_{i, j}^{(k)}$ is the length of the orthogonal projection of the convex body $P_{i, k}$ onto the vector~$v_j^*$.
In particular, 
$q_{ij}^{(0)} = \max_{z \in \cH_{i}}|(v^*_j, z)| = \max_{ l = 1, \ldots , n_i}|(v^*_j, v_i^{(l)})|$.
Sometimes we omit the superscript $k$ if its value is specified.
For given $j \in \{1, \ldots , r\}$ and for a point $x \in \re^d$, we denote
\begin{equation}\label{qjx}
q_j^{(k)}(x) \ =\
\max_{\tilde \Pi \in \tilde \cA^p, \, p = 0, \ldots, k} \bigl|\, \bigl(v^*_j\, , \, \tilde \Pi x\bigr)\, \bigr| ; \quad q_j(x) \ =\ q_j^{(\infty)}(x) \ = \
\sup_{\tilde \Pi \in \tilde \cA^p, \, p \ge 0} \bigl|\, \bigl(v^*_j\, , \, \tilde \Pi x\bigr)\, \bigr|
\end{equation}
Note that for $x = v_i$, the sets 
$\{\tilde \Pi x \ | \ \tilde \Pi \in \tilde \cA^p, \, p = 0, \ldots, k\}$
and $\cV_{i, k}$ have the same symmetrized convex hull~$P_{i, k}$. Therefore, the maxima  of 
the function $\bigl| \bigl(v^*_j\, , \, z\bigr) \bigr|$ over  these two sets coincide with
the maximum over $P_{i, k}$. Hence, comparing~(\ref{qij}) and~(\ref{qjx}) gives  
$$
 q_j^{(k)}(v_i) \ = \ q_{ij}^{(k)}\, , \ k \ge 0; \qquad 
       q_j(v_i) \ = \ q_{ij}^{(\infty)}\, . 
$$
For an arbitrary balancing  vector $\pa = (\alpha_1, \ldots , \alpha_r)$, we write 
$\pa \, \cdot \, \cH \, = \, \{\alpha_1 \cH_1, \ldots , \alpha_r \cH_r\}$. 
Our aim is to find a balancing vector such that the polytope algorithm starting simultaneously with the roots $\pa \cH$ terminates within finite time.

\begin{defi}\label{d10}
Let $k \in \{0\}\cup \n \cup \{\infty\}$. A balancing vector $\pa \in \re^r_{+}$ is called $k$-admissible, if
\begin{equation}\label{k-admiss}
\alpha_i \, q^{(k)}_{ij} \ < \ \alpha_j \ \, \qquad i, j = 1, \ldots , r, \  i \ne j\, .
\end{equation}
An $\infty$-admissible vector is called admissible.
\end{defi}
Since the value $q^{(k)}_{ij}$ is non-decreasing in $k$, we see that the $k$-admissibility for some $k$ implies the same holds true for all smaller $k$. In particular, an admissible vector is $k$-admissible for all $k \ge 0$.

We begin with two auxiliary facts needed in the proofs of our main results.

\begin{lemma}\label{l20}
If a $\, d\times d$ matrix $A$ and a vector $x \in \re^d , \, x \ne 0$, are such that
 $\|Ax - x\| < \varepsilon\, \|x\|$, then $A$ has an eigenvalue $\lambda \in \co$ for which
$|\lambda - 1| < C(d) \, \|A\|\, \varepsilon^{1/d}$, where $C(d)$ depends only on the dimension $d$.
\end{lemma}
See~\cite{SS90} for the proof. The following combinatorial fact is well-known:
\begin{lemma}\label{l30}
Let $a, b$ be two simple nonempty words of a finite alphabet, $n, k \ge 2$ be natural numbers.
If the word $a^n$ contains a subword $b^k$ such that $|b^k| > |a|$, then $b$ is a cyclic permutation of $a$.
\end{lemma}
Now we extend the key property of dominance to a set of candidate s.m.p.'s.
\begin{defi}\label{d20}
Products $\Pi_1, \ldots , \Pi_r$ are called dominant for the family $\cA$
if all numbers $[\rho(\Pi_i)]^{1/n_i}, \, i = 1, \ldots , r$, are equal (denote this value by~$\rho_c$) and
there is a constant $\gamma < 1$ such that the spectral radius of each product of matrices
from the normalized family $\tilde \cA = \rho_c^{-1}\cA$ which is neither a power of some~$\tilde \Pi_i$ nor that of its cyclic permutation is smaller than~$\gamma$.
\end{defi}

\subsection{Criterion of convergence of Algorithm~\ref{algoP}}\label{s3.2}

If the products $\Pi_1, \ldots , \Pi_r$ are dominant, then they are s.m.p., but, in general, not vice versa.
The s.m.p. property means that the function $f(\tilde \Pi)  \, = \,  [\rho(\tilde \Pi)]^{1/n}$ ($n$ is the length of~$\Pi$)
defined on the set of products of the normalized family~$\tilde \cA$
attains its maximum (equal to one) at the products $\tilde \Pi_j$ and at their powers and cyclic permutations.  The dominance property means, in addition, that for all other products, this function is smaller than some $\gamma < 1$. This property seems to be too strong, however, the following theorem shows that it is rather general.

\smallskip

\begin{theorem}\label{th10}
Algorithm~\ref{algoP} with the initial products $\Pi_1, \ldots , \Pi_r$ and with a balancing vector~$\pa$ terminates within finite time if and only if these products are all dominant, their leading eigenvalues are unique and simple and~$\pa$ is admissible.
\end{theorem}

The proof is in Appendix 1.

\begin{remark}\label{r30} {\em Theorem~\ref{th10} implies that if the algorithm
terminates within finite time, then the leading eigenvalues of products~$\Pi_i$ must be
unique and simple. That is why we defined admissible balancing  vectors for this case only.}
\end{remark}

If Algorithm~\ref{algoP} produces an invariant polytope, then our candidate s.m.p.'s are
not only s.m.p.'s but also dominant products. A number of numerical experiments suggests that the situation when the algorithm terminates within finite time (and hence, there are dominant products) should be generic.

\subsection{The existence of an admissible  balancing vector}\label{s3.3}

By Theorem~\ref{th10}, if all our candidate s.m.p.'s $\{\tilde \Pi_i\}_{i =1}^r$ are dominant and have unique simple leading eigenvalues, then balancing the
corresponding roots $\{\cH_i\}_{i =1}^r$
by weights  $\alpha_1, \ldots , \alpha_r$ we run the algorithm and construct an invariant polytope, provided the balancing vector $\alpha$ is admissible. A natural question arises if an admissible vector always exists. 
The next theorem gives an affirmative answer. Before we formulate it, we need an auxiliary result.

\begin{lemma}\label{l10}
For given coefficients $q_{ij}^{(k)}$ the system~\eqref{k-admiss} has a
solution~$\alpha > 0$ is and only if for every nontrivial cycle $(i_1, \ldots , i_n)$
on the set $\{1, \ldots , r\}$, we have (with $i_{n+1} = i_1$)
\begin{equation}\label{cycle}
\prod_{s=1}^{n} q^{(k)}_{i_{s}i_{s+1}} \ < \ 1 \, .
\end{equation}
\end{lemma}
{\tt Proof.} The necessity is simple: for an arbitrary cycle we multiply the $n$ inequalities $\alpha_{i_{s}} q^{(k)}_{i_{s}i_{s+1}} \, < \, \alpha_{i_{s+1}}, \ s = 1, \ldots , n$, and obtain~(\ref{cycle}).
To prove sufficiency, we slightly increase all numbers $q^{(k)}_{ij}$ so
that~(\ref{cycle}) still holds for all cycles.
This is possible, because the total number of cycles is finite. We set $\alpha_1 = 1$ and $\alpha_j = \max \, \prod_{s=1}^{n} q^{(k)}_{i_{s}i_{s+1}}$, where maximum is computed aver all paths $i_1 \to \cdots \to i_n \to i_{n+1}$ with $i_1 = 1, i_{n+1} = j, \, n \ge 0$.
Note that if a path contains a cycle, then removing it increases the product
$\prod_{s=1}^{n} q^{(k)}_{i_{s}i_{s+1}}$, since the corresponding product along the cycle
is smaller than one. This means that, in the definition of $\alpha_j$,
it suffices to take the maximum over all simple (without repeated vertices) paths, i.e., over a finite set.

It is easy to see that $\alpha_{i} q^{(k)}_{ij} \, \le \, \alpha_{j}$.
Reducing now all $q^{(k)}_{ij}$ back to the original values, we obtain strict inequalities.
{\hfill $\Box$}
\smallskip

\begin{theorem}\label{th20}
If  the products $\Pi_1, \ldots \Pi_r$  are dominant and have unique simple leading eigenvalues, then they have an admissible balancing vector.
\end{theorem}
{\tt Proof.} In view of Lemma~\ref{l10},  it suffices to show that for every cycle
$(i_1, \ldots,  i_n), \, n \ge 2$,  on the set $\{1, \ldots , r\}$, we have
$\prod_{s=1}^{n} q_{i_{s}i_{s+1}} \ < \ 1 $. We denote this quantity by $h$ and take arbitrary 
$\delta > 0$. There is a product $\tilde \Pi$ of matrices from $\tilde \cA$ such that 
$|\bigl( v_{i_2}^*, \tilde \Pi v_{i_1}\bigr)| >  q_{i_{1}i_{2}} - \delta$.
Without loss of generality we assume that this scalar product is positive.
Since the product $\tilde \Pi_{i_2}$ has a unique simple leading  eigenvalue $1$, it follows that for every $x\in \re^d$ we have $\tilde \Pi_{i_2}^k x \to (v_{i_2}^*, x)\, v_{i_2}$ as $k \to \infty$.
Applying this to the vector $x = \tilde \Pi v_{i_1}$, we conclude that
$\|\tilde \Pi_{i_2}^k \tilde \Pi v_{i_1} - q_{i_{1}i_2} v_{i_2}\| < 2\delta$, whenever $k$ is large enough. Thus, for the product $\tilde \Pi^{(1)} = \tilde \Pi_{i_2}^k \tilde \Pi$, the vector 
$\tilde \Pi^{(1)}v_{i_1}$ is close to $q_{i_{1}i_{2}} v_{i_2}$. 
Analogously, for each $s = 1, \ldots , n$, we find a
product $\tilde \Pi^{(s)}$ such that  the vector $\tilde \Pi^{(s)}v_{i_s}$ is close
to $q_{i_{s}i_{s+1}} v_{i_{s+1}}$. Therefore, for the product
$B  = \prod_{s=1}^n \tilde \Pi^{(s)}$, the vector $Bv_{i_1}$ is close to
$\bigl( \prod_{s=1}^n q_{i_{s}i_{s+1}}\bigr) \, v_{i_1} \, = \, h \, v_{i_1}$.
Note that $\|B\| \le M$, where~$M$ is the supremum of norms of all products of matrices from~$\tilde A$. 
If $h \ge 1$, then invoking Lemma~\ref{l20} we conclude that
$\rho(B) \ge 1 - \varepsilon$, where $\varepsilon > 0$ can be made 
arbitrarily small by taking $k \to \infty$.  Due to the dominance assumption, it follows that
$B$ is a power of some $\tilde \Pi_0 \in \Omega$, where $\Omega$ is the set of products 
$\tilde \Pi_1, \ldots , \tilde \Pi_r$ and of its cyclic permutations. 
Due to the dominance assumption, it follows that $B$ is a power of some
$\tilde \Pi_0 \in \Omega$.
Taking $k$ large enough we apply Lemma~\ref{l30} to the words
$a = \tilde \Pi_0, b = \tilde \Pi_{i_2}$ and conclude that $\tilde \Pi_{i_2}$ is a cyclic permutation of $\tilde \Pi_0$. Similarly, $\tilde \Pi_{i_3}$ is a cyclic permutation of $\tilde \Pi_0$. This is impossible, because $i_2 \ne i_3$, and the products $\tilde \Pi_{i_2}, \tilde \Pi_{i_3}$ are not cyclic permutations of each other. The contradiction proves that $h < 1$ which completes the proof of the theorem.
{\hfill $\Box$}
\smallskip

\begin{remark}\label{r100} 
{\em  In a just published paper~\cite{MR14}, M\"oller and Reif present another approach 
for the computation of joint spectral radius. Developing ideas from~\cite{HMR} they come up with an elegant branch-and-bound algorithm, which, in contrast to the classical 
branch-and-bound method~\cite{G}, can find the exact value. Although its running time is typically bigger than for our invariant polytope algorithm~\cite{GP13}
(we compare two algorithms in Example \eqref{ex:8p} below), it has several advantages. In particular, it uses the same scheme for the cases of one and of many s.m.p. It would be interesting to analyze possible application of the balancing idea 
for that algorithm.}
\end{remark}

\subsection{How to find the balancing vector}\label{s3.4}

Thus, an admissible balancing vector $\pa$ does exist, provided our candidate s.m.p.'s are dominant products and their leading eigenvectors are unique and simple.
To find $\pa$, we take some $k \ge 0$, compute the values $q^{(k)}_{ij}$ by evaluating polytopes $P_{i, k}, \, i = 1, \ldots , r$, set
$y_i = \log \alpha_i\, , \ b_{ij}^{(k)} = - \log q^{(k)}_{ij}$ and solve the following LP problem with variables $y_0, \ldots , y_r$:
\begin{equation}\label{LP}
\left\{
\begin{array}{l}
y_0 \ \to \ \max\\
y_i\, - \, y_j \ \le \, - y_0 \, + \, b^{(k)}_{ij}, \qquad i, j = 1, \ldots , r, \  i \ne j\, .
\end{array}
\right.
\end{equation}
If $y_0 \le 0$, then the $k$-admissible vector does not exist. In this case,  we have to either increase $k$ or find other candidate s.m.p.'s. If $y_0 > 0$, then we have a $k$-admissible vector $\pa = (e^{y_1}, \ldots , e^{y_r})$. This vector is optimal in a sense
that the minimal ratio between $\frac{\alpha_j}{\alpha_i}$ and $q^{(k)}_{ij}$ over all $i, j$ is the biggest possible.

\begin{remark}\label{r10}
{\em To find an admissible vector one needs to solve LP problem~(\ref{LP}) for $k = \infty$.
However, in this case the evaluation of the coefficients $b^{(k)}_{ij}$ may, a priori, require an infinite time.
Therefore, we solve this problem for some finite $k$ and then run Algorithm~\ref{algoP} with the
obtained balancing vector~$\pa$. If the algorithm terminates within finite time, then $\pa$ is
admissible indeed (Theorem~\ref{th10}). Otherwise, we cannot conclude that there are no admissible balancing and that our candidate s.m.p.'s are not dominant. We try to to increase $k$ and find a new vector~$\pa$.
}
\end{remark}

Thus, {\em Step 5} of Algorithm~\ref{algoP} consists in choosing a reasonably big~$k$ and
solving LP problem~(\ref{LP}). If it results $y_0 \le 0$, then the balancing vector does not exist, and hence the algorithm will never converge and we have to find another candidate s.m.p.  If $y_0 > 0$, then the vector $\pa = (e^{y_1}, \ldots , e^{y_r})$ is $k$-admissible. If the algorithm does not converge with this $\alpha$, we increase~$k$ and solve~(\ref{LP}) again.
\smallskip

\begin{remark}\label{r15}
{\em Our approach works well also if the family has a unique s.m.p.~$\Pi_1$, but
there are other simple products~$\Pi_2, \ldots , \Pi_r$
for which the values $[\rho(\tilde \Pi_j)]^{1/n_j}$, although being smaller than
$[\rho(\tilde \Pi_1)]^{1/n_1} = 1$, are close to it.
In this case the (original) invariant polytope algorithm sometimes converges slowly performing many iterations and producing many vertices. This is natural, because if,
say $[\rho(\tilde \Pi_2)]^{1/n_2} = 1-\delta$ with very small $\delta$>0, then the dominance of 
$\tilde \Pi_1$ over $\tilde \Pi_2$ plays a role only after many iterations. Our approach suggests
to collect all those ``almost s.m.p. candidates''~$\Pi_2, \ldots , \Pi_r$ add them
to $\Pi_1$, find the balancing multipliers $\{\alpha_i\}_{i=1}^r$ for their roots
by solving LP problem~(\ref{LP}) and run Algorithm~\ref{algoP} for the initial
set $\cV_0 = \{\alpha_j\cH_j\}_{j=1}^r$.  In most of practical cases, this modification significantly speeds up the algorithm.
\smallskip

Another modification of the invariant polytope algorithm is considered in the next section.
}
\end{remark}

\begin{ex}
\label{ex:8p}

\rm Consider the following example introduced by Deslaurier and Dubuc in \cite{DD89}, associated to an eight-point subdivision scheme,
\begin{eqnarray*}
A_1 & = & \left( \begin{array}{rrrrrrrr}
30  & -14 &  -14  & 30  &  0  &  0  &  0  &  0 \\
-5  & -56 &  154 &  -56  & -5  &  0  &  0  &  0 \\
 0  & 30  & -14  & -14  & 30  &  0  &  0  &  0 \\
 0  & -5  & -56  & 154  &  -56  & -5  &  0  &  0 \\
 0  &  0  & 30  & -14  & -14  & 30  &  0  &  0 \\
 0  &  0  & -5  & -56 &  154 &  -56  & -5  &  0 \\
 0  &  0  &  0  & 30  & -14  & -14  & 30  &  0 \\
 0  &  0  &  0  & -5  & -56 &  154  & -56  & -5
\end{array} \right), 
\\
A_2 & = & \left( \begin{array}{rrrrrrrr}
-5  & -56 &  154 &  -56  & -5  &  0  &  0  &  0 \\
 0  & 30  & -14  & -14  & 30  &  0  &  0  &  0 \\
 0  & -5  & -56  & 154 &  -56  & -5  &  0  &  0 \\
 0  &  0  & 30  & -14  & -14  & 30  &  0  &  0 \\
 0  &  0  & -5  & -56  & 154 &  -56  & -5  &  0 \\
 0  &  0  &  0  & 30  & -14  & -14  & 30  &  0 \\
 0  &  0  &  0  & -5  & -56  & 154  & -56  & -5 \\
 0  &  0  &  0  &  0  &  30  & -14  & -14  & 30
\end{array} \right).
\end{eqnarray*}

The joint spectral radius of $\cA = \{A_1,A_2\}$ was found in \cite{MR14}, where 
it was shown that both $A_1$ and $A_2$ are s.m.p. 
Its computation required the construction of a binary tree with $14$ levels and considering about $130$ matrix products (i.e., vertices of the tree). Applying our Algorithm~\ref{algoP} with the candidates s.m.p.  $A_1$ and $A_2$ and with a balancing vector $\alpha = ( 1 \quad 1)^{\rm T}$ for the leading eigenvectors $v_1$ and $v_2$ of
$A_1$ and $A_2$, respectively, we construct the invariant polytope with $24$ vertices 
in $5$ steps:
\[
\begin{array}{llllll}
v_1  &  v_2  &  v_3 = \tA_1\,v_2  &  
v_4 = \tA_2\,v_1  &  v_5 = \tA_1\,v_3 & v_6 = \tA_1\,v_4 \\
v_7 = \tA_2\,v_3  &  v_8 = \tA_2\,v_4  &  v_9 = \tA_1\,v_5  &  
v_{10} = \tA_1\,v_6  &  v_{11} = \tA_1\,v_7 & v_{12} = \tA_1\,v_8 \\
v_{13} = \tA_2\,v_5  &  v_{14} = \tA_2\,v_6  &  v_{15} = \tA_2\,v_7  &  
v_{16} = \tA_2\,v_8  &  v_{17} = \tA_1\,v_{11} & v_{18} = \tA_1\,v_{12} \\
v_{19} = \tA_1\,v_{13}  &  v_{20} = \tA_1\,v_{14}  &  v_{21} = \tA_2\,v_{11}  &  
v_{22} = \tA_2\,v_{12}  &  v_{23} = \tA_2\,v_{13} & v_{24} = \tA_2\,v_{14}. 
\end{array}
\]
Thus, in our case it suffices to construct a binary tree with $5$ levels and consider $24$ of its vertices.
\end{ex}
\smallskip

\section{Introducing extra initial vertices}\label{s4}
\smallskip

The same approach developed for the case of many s.m.p. can be used to
introduce extra initial vertices. Sometimes Algorithm~\ref{algoP} converges slowly
because the family $\cA$ is not well-conditioned: its matrices have a common ``almost
invariant subspace'' of some dimension~$s \le d-1$.   In this case the invariant polytope~$P$ may be very flattened (almost contained in that subspace). As a consequence, the algorithm performs many iterations because the polytopes~$P_k$, being all flattened, badly  absorb new vertices. To avoid this trouble one can introduce extra
initial vertices $x_1, \ldots , x_s$ and run Algorithm~\ref{algoP} with the initial set
$\cV_0 = \{\alpha_1 \cH_1, \ldots , \alpha_r \cH_r, x_1, \ldots , x_s\}$.
In next theorem we use the value~$q_j(x)$ defined in~(\ref{qjx}).

\begin{theorem}\label{th30}
Suppose Algorithm~\ref{algoP} with initial roots $\cH_1, \ldots , \cH_r$ terminates within finite time; then this algorithm with extra initial vertices $x_1, \ldots , x_s$ also does if and only if $q_j(x_i) < 1$ for all $j = 1, \ldots , r; \, i = 1, \ldots , s$.
\end{theorem}
{\tt Proof.} For the sake of simplicity, we consider the case $s=1, r =1$, the proof in the general case is similar.
Let $P$ denote the final polytope produced by the algorithm starting with the
root $\cH_1$, and let
$\, P(x_1) \, = \, {\rm absco} \,
\{\tilde \Pi x_1 \ | \ \tilde \Pi \in \tilde \cA^n, \, n \ge 0\}$.

Necessity. Assume the algorithm terminates within finite time.
In the proof of Theorem~\ref{th10} we showed that the maximum of the linear functional
$f(x) = (v_1^*, x)$ on the final polytope~$P$ is equal to one and is attained at a unique point $v_1$.
  Hence, either $q_1(x_1) <1$, in which case the proof is completed, or
 $q_1(x_1) = 1$,  and hence there is a sequence of products $\{\tilde \Pi^{(k)}\}_{k \in \n}$ such that $(v_1^*, \tilde \Pi^{(k)}x_1) \to 1$
 as $k \to \infty$. Therefore, $\tilde \Pi^{(k)}x_1 \to v_1$ as $k \to \infty$.
 This implies that, for sufficiently large $k$, the points $\tilde \Pi^{(k)}x_1$ are not absorbed in the algorithm, and hence, the algorithm cannot terminate within finite time.

 Sufficiency. Since the second largest eigenvalue of the matrix~$\tilde \Pi_1$ is smaller than~$1$ in absolute value, it follows that
 $\tilde \Pi_1^k \, \to \, v_1\, [v_1^*]^T$ as $k \to \infty$. If $q_1(x_1) < 1$, then
 the matrix $v_1 \, [v_1^*]^T$ maps the set  $P(x_1)$ to the segment 
$[-q_1(x_1)\, v_1 \, , \, q_1(x_1)\, v_1]$, which is contained in $q_1(x_1)\, P$. 
Hence, $\tilde \Pi_1^k\, (P(x_1)) \, \subset \, P$, for some~$k$. Therefore, every product $\tilde \Pi$ of length~$N$ containing a subword $\tilde \Pi_1^k$ takes the point $x_1$ inside $P$. On the other hand, for all products $\tilde \Pi$ not containing this subword, we have $\|\tilde \Pi\| \le C q^{N}$, where $C > 0, \, q \in (0,1)$ are some constants (see~\cite[Theorem~4]{GP13}).
  Hence, for large $N$, all such products also take the point~$x_1$ inside~$P$. Thus, all long products take $x_1$ inside $P$, hence the algorithm starting with the initial set $\cH_1 \cup \{x_1\}$ terminates within finite time.
{\hfill $\Box$}
\smallskip

In practice,  it suffices to introduce extra initial vertices of the form
$x_i = \mu_i e_i$, where $e_i$ is the canonical basis vector (the $i$-th entry is one and all others are zeros) and $\mu_i$ is some positive coefficient. We fix a reasonably small $\varepsilon > 0$, say, between  $0.001$ and $0.1$, a reasonably large $k$, 
say $k=15$, rum $k$ iterations of Algorithm~\ref{algoP} and compute the values
$$
Q_i^{(k)} \ = \
\max_{v \in \cV_k} \ \bigl|\, \bigl(e_i\, , \, v\bigr)\, \bigr|\ , \quad i = 1, \ldots , d\, .
$$
So, $Q_i^{(k)}$ is the length of projection of the polytope~$P_k$ onto the $i$-th coordinate axis, or the largest $i$-th coordinate of its vertices.
If $Q_i^{(k)} \ge \varepsilon$ for all $i$, then $P_k$ contains the cross-polytope
${\rm absco} \, \{ \,  \varepsilon e_i\ ,  \ i = 1, \ldots , d\}$, which, in turn, contains the Euclidean ball of radius $\varepsilon/\sqrt{d}$ centered at the origin. In this case the polytope $P_k$ is considered to be well-conditioned, and hence we do not add any extra vertex.
If, otherwise, $Q_i^{(k)} < \varepsilon$ for some $i$, then we add an extra vertex
$x_i = \varepsilon \, [\max_j q^{(k)}_j(e_i)]^{-1}e_i$. Collecting all such  vertices for $i = 1, \ldots , d$
(assume there are $s \le d$ ones) we run Algorithm~\ref{algoP} with $s$ initial extra vertices in the set~$\cV_0$.

In Section~\ref{s6} we apply  this trick to speed up Algorithm~\ref{algoP} for Daubechies matrices, which turn out to be extremely ill-conditioned.

\section{Applications:  the Butterfly subdivision scheme}\label{s5}

Subdivision schemes are iterative algorithms of linear interpolation and approximation
of multivariate functions and of generating curves and surfaces. Due to their remarkable properties they are widely implemented and studied in an extensive literature.

The Butterfly scheme originated with Dyn, Gregory, and Levin~\cite{DGL90} and became one of the most popular bivariate schemes for  interpolation and for generating smooth surfaces (see also the generalization given in \cite{ZZ}).
This scheme is a generalization of the univariate four-point interpolatory scheme
 to bivariate functions~\cite{AD, HMR, SDL}. First we take an arbitrary triangulation of the approximated surface and consider the corresponding piecewise-linear interpolation.
 This interpolation produces a sequence of piecewise-linear surfaces with thriangular faces that converges to a continuous surface, which is  considered as an interpolation of the original one.
 To describe this algorithm in more detail we assume that the original surface is given by a bivariate function $f(x_1, x_2)$. We consider a regular triangulation of $\re^2$
 and take the values of the function~$f$ at its vertices. So we obtain a function $f_1$
 defined on a triangular mesh. In the next iteration we define the function $f_2$ on the
 refined triangular mesh: the values at the vertices of the original mesh stay the same,
 the values at midpoints of edges are defined as a linear combination
 of eighth neighboring vertices as shown in the following figure ($X$ is the new vertex,
 the coefficients of the linear combination are written in the corresponding vertices).
\begin{figure}[ht]
\centering \global\def\path{#1}\input{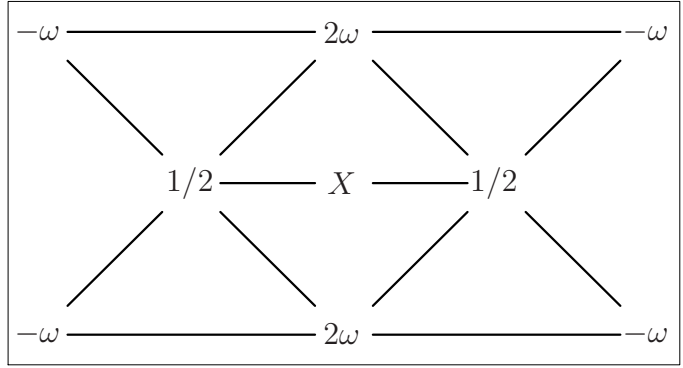}
\caption{The generalized Butterfly scheme. \label{fig:GBS}}
\end{figure}
 The parameter $\omega$ is the same for all vertices and for all iterations.
 In the next iteration we do the same with the new mesh and produce the function $f_3$, etc.
 Each function is extended from the corresponding mesh to the whole plane by linearity.
The scheme is said to converge if for every initial function~$f$, the functions $f_k$
converge uniformly to a continuous function~$\cS(f)$. Due to linearity and shift-invariance,
it suffices to have the convergence  for the initial
$\delta$-function, which vanishes at all vertices but one, where it is equal to one.
The corresponding limit function $\cS(\delta)$ is called  {\em refinable function} and denoted by 
$\varphi$.
The scheme converges for every $\omega \in \bigl(0, \frac14 \bigr)$ and reproduces polynomial of degree one, i.e., if $f(x_1, x_2) = a_1x_1 + a_2x_2 + a_0$, then
$\cS(f) = f$. The case $\omega = 1/16$ is special, in this case the scheme reproduces polynomials of degree~$3$~\cite{DGL90}.
One of the most important problems in the study of any subdivision scheme is its regularity.
We use the standard modulus of continuity
$$
\omega_{f}(h) \ = \ \sup\, \bigl\{ |f(x+\xi) - f(x)| \quad | \quad  x \in \re^d\, , \|\xi\| \le h \bigr\}\, .
$$
The {\em H\"older exponent} of the function~$f$ is
$$
\alpha_{f}(h) \ = \ n \ + \ \sup \bigl\{ \alpha \ge 0 \ \bigl|  \
\omega_{f^{(n)}} (h) \, \le \, C h^{\, \alpha} \, , \ h > 0\,  \bigr\}\, ,
$$
where $n$ is the biggest integer such that $f \in C^n(\re^d)$.
The H\"older regularity of a subdivision scheme is $\alpha_{\varphi}$.
This is well-known that the exponent of H\"older regularity of a 
bivariate subdivision scheme is equal to
$$
\alpha_{\varphi} \ = \ - \,\log_2 \rho \bigl( T_1^{(\ell)}, T_2^{(\ell)},
T_3^{(\ell)}, T_4^{(\ell)} \bigr),
$$
where $\ell$ is the maximal degree of the space of algebraic polynomials reproduced by the scheme, $T_i^{(\ell)}$ are restrictions of the {\em transition operators} $T_i$ of the scheme to their common invariant  subspace orthogonal to the subspace of algebraic polynomials of degree~$\ell$~\cite{NPS}. For the Butterfly scheme, for all $\omega \ne 1/16$, we have $\ell = 1$ and the operators $T_i^{(1)}$ are given by $24 \times 24$-matrices. In the only ``most regular'' case $\omega = 1/16$ the scheme respects the cubic polynomials, i.e., $\ell=3$, and $T_i^{(3)}$ are $17 \times 17$-matrices.
For this exceptional and most important case it was conjectured in early 90th that
$\alpha_{\varphi} = 2$,  i.e., the limit function~$\varphi$ of the scheme are continuously differentiable and its derivative~$\varphi\,'$ has H\"older exponent $1$. We are going to prove this conjecture and, moreover, we show that the derivative of the function $\varphi$ is not Lipschitz, but ``almost Lipschitz'' with the logarithmic
factor~$2$.
\begin{theorem}\label{th40}
The H\"older regularity of the Butterfly scheme with $\omega = \frac{1}{16}$ is equal
to~$2$.
The derivative $\varphi\,'$ of the limit function is ``almost Lipschitz'' with the logarithmic 
factor $2$:
\begin{equation}\label{almostL}
\omega_{\varphi'}(h) \, \asymp \, h \, |\log  h |^2\, , \quad  h \in \Bigl(0, \frac12 \Bigr)\, .
\end{equation}
\end{theorem}
\begin{remark}\label{r20}
{\em In the four-point subdivision scheme, which is a univariate parameter-dependent analogue of the butterfly scheme, the case $\omega = \frac{1}{16}$ is also crucial.
This is the only case when the scheme reproduces cubic polynomials. As it was proved by S.Dubuc 
in 1986~\cite{D86}, the regularity of the four-point scheme in this case is equal to two and 
$\omega_{\varphi\,'} (h) \asymp h |\log h|$, i.e., $\varphi\,'$ is almost Lipschitz with the logarithmic factor~$1$.
By Theorem~\ref{th40}, for the Butterfly scheme, the situation is similar, but
$\varphi\,'$ is almost Lipschitz with the logarithmic factor~$2$.}
\end{remark}
To prove Theorem~\ref{th40} we first show that $\rho \{T_i^{(3)}, \ i = 1, \ldots , 4\} = \frac14$. Then we conclude that
$\rho \{T_i^{(1)}, \ i = 1, \ldots , 4\} = \frac14$. By a more refine analysis of the matrices we establish that
\begin{equation}\label{defect1}
\max\, \Bigl\{   \| T_{i_k}^{(1)} \cdots  T_{i_1}^{(1)}\|  \ \Bigl| \
i_1, \ldots , i_k \in \{1,2,3,4\}\,   \Bigr\} \ \asymp \ k^2 \, 4^{-k}\, .
\end{equation}
Then it will remain to refer to some known facts of the theory of subdivision schemes.
The main and most difficult part is the finding of the joint spectral radius of the matrices 
$T_i^{(3)}$. This is done in the next subsections. Then we conclude the proof.

\subsection{The case $\omega=1/16$: the classical Butterfly scheme}\label{s5.1}

The $17 \times 17$-matrices $T_i^{(3)}$ can be computed exactly by the Matlab program 
of P.Oswald \cite{O12}. 
To simplify the notation, we denote $A_i = 4 T_i^{(3)}\, , \, i = 1, \ldots , 4$,
and $\cA = \{A_1, \ldots , A_4\}$. These four $17 \times 17$-matrices  are written in Appendix 2.
Our goal is to prove that  $\rho (\cA) = 1$.
It is remarkable that all the matrices of the family~$\cA$ and the leading eigenvector of its s.m.p. possess rational entries, so our computations are
actually done in the exact arithmetics.
\smallskip

\textbf{Step 1. Factorization of the family~$\cA$ to $\cA_1$ and $\cA_2$}.
\smallskip

It appears that the matrices $A_1, \ldots , A_4$ can be factored to a block lower-triangular 
form in a common basis.
To see this we take the leading eigenvectors $v_1$ of the matrix $A_1$
(as it was mentioned above, it has rational entries):
$$
\begin{array}{l}
v_1 \ = \ \Bigl(
 -\frac{1746890}{2004757}, \,
 -\frac{942260}{2004757}, \,
 -\frac{1391945}{2004757}, \,
 -\frac{1596995}{2004757}, \,
 -\frac{1987045}{4009514}, \,
 -\frac{3186205}{4009514}, \,
 \frac{1392363}{4009514}, \, \\
 \frac{1478789}{2004757}, \,
 \frac{2902096}{2004757}, \,
 \frac{594645}{4009514}, \,
 \frac{1063787}{4009514}, \,
 \frac{2802569}{4009514}, \,
 \frac{2242099}{4009514}, \,
 -\frac{45664}{2004757}, \,
 0, \,
 0, \,
 1
\Bigr)^T\, .
\end{array}
$$
It is checked directly that the linear span of the following six vectors  is a common invariant subspace of all matrices from~$\cA$:
$v_1, \, v_2 = A_2 v_1, \, v_3 = A_3 v_1, \, v_4 = A_4 v_1, \, v_5 = A_1 A_2 v_1, \, v_6 = A_1 A_3 v_1$.
Therefore,  we can transform the matrices from $\cA$ into block lower-triangular form
with diagonal blocks of dimensions $6$ and $11$.
 The transformation matrix is
\[
S = \left( \begin{array}{cccccccccc}
e_1 & e_2 & \ldots & e_{11} & v_1 & v_2 & v_3 & v_4 & v_5 & v_6 \\
\end{array}
\right)
\]
(written by columns), where $e_k \in \R^{17}$ is the $k$-th canonical basis vector. This gives the transformed matrices with block lower-triangular structure:
\[
S^{-1} A_i S = \left( \begin{array}{rr} C_i & 0 \\ D_i & B_i \end{array} \right), \qquad i=1,\ldots,4.
\] \label{BC}
with $6\times 6$-matrices $C_i$ and $11\times 11$-matrices $B_i$. Those matrices are written  down in Appendix 2. Note that they are all rational.
Let $\cA_1 = \{ B_1, B_2, B_3, B_4 \}$ and $\cA_2 = \{ C_1, C_2, C_3, C_4 \}$.
It is well known that the joint spectral radius of a block lower-triangular family of matrices is equal to the maximal joint spectral radius of blocks~\cite{BW92}. 
Hence $\, \rho(\cA)  = \max \, \{\rho(\cA_1), \rho(\cA_2)\}$.
Now we are going to show that $\rho(\cA_1) = \rho(\cA_2) = 1$, from which it will follow that
$\rho(\cA) = 1$. We begin with the family~$\cA_1$.
\smallskip

\textbf{Step 2.
Analysis of the family $\cA_1$}

\smallskip

We have a family $\cA_1$ of $11\times 11$-matrices $B_1, \ldots , B_4$ written in Appendix 2.
Each of the  matrices~$B_1, B_2, B_3$ has a simple leading eigenvalue $1$,
the corresponding leading eigenvectors $u_1, u_2, u_3$ are all simple. The  matrix
 $B_4$ has spectral radius
$1/2$. We are going to show that $\rho (\cA_1) = 1$, i.e., this family has three s.m.p.: $B_1$, $B_2$ and $B_3$.
The leading eigenvectors are (normalized in the maximum norm),
\[
u_1 = \left( \begin{array}{r}
 \frac{2497}{3306} \\[2mm]
 \frac{1453}{3306}  \\[2mm]
-\frac{1997}{3306}  \\[2mm]
 \frac{283}{3306}  \\[2mm]
 \frac{2023}{3306} \\[2mm]
 \frac{787}{3306}  \\[2mm]
 \frac{2851}{11020} \\[2mm]
-\frac{4749}{11020} \\[2mm]
-\frac{613}{2204} \\[2mm]
-1 \\[2mm]
1
\end{array} \right), \qquad
u_2 = \left( \begin{array}{r}
 \frac{230}{1089} \\[2mm]
 \frac{230}{1089} \\[2mm]
-\frac{370}{1089} \\[2mm]
-\frac{730}{1089} \\[2mm]
-\frac{730}{1089} \\[2mm]
-\frac{370}{1089} \\[2mm]
-\frac{221}{363}  \\[2mm]
 1 \\[2mm]
-\frac{221}{363}  \\[2mm]
 0 \\
 0
\end{array} \right), \qquad
u_3 = \left( \begin{array}{r}
 \frac{66}{625} \\[2mm]
 \frac{1158}{4375} \\[2mm]
 \frac{2298}{4375} \\[2mm]
 \frac{82}{4375} \\[2mm]
-\frac{154}{625} \\[2mm]
 \frac{442}{4375} \\[2mm]
-1 \\[2mm]
 \frac{2661}{21875} \\[2mm]
-\frac{15959}{21875} \\[2mm]
 \frac{2204}{4375} \\[2mm]
-\frac{2204}{4375}
\end{array} \right).
\]
Solving problem~(\ref{LP}) for $r=3, k=10$, and for the family~$\cA_1$, we obtain the scaling factors:
\[
\alpha_1 = 0.50379\ldots, \qquad  \alpha_2 =  0.48126\ldots, \qquad \alpha_3 = 1
\]
The Algorithm~\ref{algoP} with the three candidate s.m.p.'s
$\Pi_1 = B_1, \Pi_2 = B_2$, and $\Pi_3 = B_3$ and with those factors~$\alpha_i$ terminates within four iterations and produces
an invariant polytope. However, we slightly change the factors to $\alpha_1 = 0.5$, $\alpha_2 = 0.5$ and $\alpha_3 = 1$ respectively in order to preserve the rationality of the vectors (and consequently the exactness of the computation). Thus,
\[
v_1 = \frac12 u_1, \qquad v_2 = \frac12 u_2, \qquad v_3 = u_3.
\]
The algorithm still converges within four iterations producing the polytope $P_1$ with $75\cdot 2$ vertices. Here is the list of vertices: \begin{eqnarray*}
\begin{array}{llllll}
v_1  &  v_2  &  v_3  &  v_4 = B_1\,v_2  &  v_5 = B_1\,v_3 \\
v_6 = B_2\,v_1  &  v_7 = B_2\,v_3  &  v_8 = B_3\,v_1  &  v_9 = B_3\,v_2  &  v_{10} = B_4\,v_1 \\
v_{11} = B_4\,v_2  &  v_{12} = B_4\,v_3  &  v_{13} =  B_1\,v_4  &  v_{14} = B_1\,v_5  &  v_{15} = B_1\,v_6 \\
v_{16} = B_1\,v_7  &  v_{17} = B_1\,v_8  &  v_{18} = B_1\,v_9  &  v_{19} = B_1\,v_{10}  &  v_{20} = B_1\,v_{11} \\
v_{21} = B_1\,v_{12}  &  v_{22} = B_2\,v_4  &  v_{23} = B_2\,v_5  &  v_{24} = B_2\,v_6  &  v_{25} = B_2\,v_7 \\
v_{26} = B_2\,v_8  &  v_{27} = B_2\,v_9  &  v_{28} = B_2\,v_{10}
 &  v_{29} = B_2\,v_{11}  &  v_{30} = B_2\,v_{12} \\
v_{31} = B_3\,v_4  &  v_{32} = B_3\,v_5  &  v_{33} = B_3\,v_6  &  v_{34} = B_3\,v_7  &  v_{35} = B_3\,v_8 \\
v_{36} = B_3\,v_9  &  v_{37} = B_3\,v_{10}  &  v_{38} = B_3\,v_{11}
 &  v_{39} = B_3\,v_{12}  &  v_{40} = B_4\,v_4 \\
v_{41} = B_4\,v_5  &  v_{42} = B_4\,v_6  &  v_{43} = B_4\,v_7   &  v_{44} = B_4\,v_8  &  v_{45} = B_4\,v_9 \\
v_{46} = B_4\,v_{10}  &  v_{47} = B_4\,v_{11}  &  v_{48} = B_4\,v_{12}
 &  v_{49} = B_1\,v_{15}  &  v_{50} = B_1\,v_{17} \\
v_{51} = B_1\,v_{20}  &  v_{52} = B_1\,v_{21}  &  v_{53} = B_1\,v_{30}
 &  v_{54} = B_1\,v_{38}  &  v_{55} = B_1\,v_{46} \\
v_{56} = B_2\,v_{21}  &  v_{57} = B_2\,v_{22}  &  v_{58} = B_2\,v_{27}
 &  v_{59} = B_2\,v_{28}  &  v_{60} = B_2\,v_{30} \\
v_{61} = B_2\,v_{37}  &  v_{62} = B_2\,v_{47}  &  v_{63} = B_3\,v_{20}
 &  v_{64} = B_3\,v_{28}  &  v_{65} = B_3\,v_{32} \\
v_{66} = B_3\,v_{34}  &  v_{67} = B_3\,v_{37}  &  v_{68} = B_3\,v_{38}
 &  v_{69} = B_3\,v_{48}  &  v_{70} = B_4\,v_{20} \\
v_{71} = B_4\,v_{21}  &  v_{72} = B_4\,v_{28}  &  v_{73} = B_4\,v_{30}
 &  v_{74} = B_4\,v_{37}  &  v_{75} = B_4\,v_{38}.
\end{array}
\end{eqnarray*}
Thus, $\rho(\cA_1) = 1$ and, by Theorem~\ref{th10},  $B_1, B_2, B_3$ are dominant products for $\cA_1$.

\smallskip

\textbf{Step 3. Analysis of the family $\cA_2$}
\smallskip

We have a family $\cA_2$ of $6\times 6$-matrices $C_1, \ldots , C_4$ written in
Appendix 2.
Each of the  matrices~$C_1, C_2, C_3$ has a simple leading eigenvalue $1$,
the  matrix~$C_4$ has spectral radius $1/2$.

First of all, we observe the existence of three invariant $2$-dimensional
subspaces of all the matrices~$C_i$.
We indicate by $w_1, w_2$ and $w_3$ the unique leading eigenvectors
associated to
the eigenvalue $1$ of the matrices $C_1 C_2$, $C_1 C_3$ and $C_2 C_3$,
(normalized in maximum norm),
\[
w_1 = \left( \begin{array}{r}
\frac{1}{4} \\[2mm]
\frac{1}{4} \\[2mm]
0 \\[2mm]
0 \\[2mm]
1 \\[2mm]

\end{array} \right), \qquad
w_2 = \left( \begin{array}{r}
\frac{1}{4} \\[2mm]
0 \\[2mm]
\frac{1}{4} \\[2mm]
0 \\[2mm]
0 \\[2mm]
1
\end{array} \right), \qquad
w_3 = \left( \begin{array}{r}
  -\frac{1}{7}   \\[2mm]
  -\frac{25}{28} \\[2mm]
  -\frac{25}{28} \\[2mm]
  -\frac{2}{7}   \\[2mm]
    1   \\[2mm]
    1   \\[2mm]
\end{array} \right)
\]
The invariant subspaces are given by $V_1 = {\rm span} \left( w_1, C_4
w_1 \right)$,
$V_2 = {\rm span} \left( w_2, C_4 w_2 \right)$ and
$V_3 = {\rm span} \left( w_3, C_1 w_3 \right)$. Thus we define the matrix
\[
S = \left( w_1, C_4 w_1, w_2, C_4 w_2, w_3, C_1 w_3 \right)
\]
which block-diagonalizes all matrices $C_i$, $i=1,\ldots,4$.
We denote the diagonal blocks of the matrices $S^{-1} C_i S$ as
$G_{i1}, G_{i2}$ and $G_{i3}$. We obtain three families of $2 \times 2$
matrices to analyze,
$G_1 = \{ G_{11}, G_{12}, G_{13}, G_{14} \}$ with
\[
G_{11} = \left( \begin{array}{rr}
1   & -\frac{1}{4}   \\[2mm]       0  & -\frac{1}{4}
\end{array}
\right), \qquad
G_{12} = G_{11}, \qquad
G_{13} = \left( \begin{array}{rr}
-\frac{1}{4} & 0 \\[2mm]  -\frac{1}{4} & 1
\end{array}
\right), \qquad
G_{14} = \left( \begin{array}{rr}
0 & -\frac{1}{4} \\[2mm] 1 & -\frac{1}{4}
\end{array}
\right).
\]
then $G_2 = \{ G_{21}, G_{22}, G_{23}, G_{24} \}$ with
\[
G_{21} = \left( \begin{array}{rr}
1   & -\frac{1}{4}   \\[2mm]       0  & -\frac{1}{4}
\end{array}
\right), \qquad
G_{22} = \left( \begin{array}{rr}
-\frac{1}{4} & 0 \\[2mm]  -\frac{1}{4} & 1
\end{array}
\right), \qquad
G_{23} = G_{12}, \qquad
G_{24} = \left( \begin{array}{rr}
0 & -\frac{1}{4} \\[2mm] 1 & -\frac{1}{4}
\end{array}
\right).
\]
i.e. $G_2 = G_1$ and
$G_3 = \{ G_{31}, G_{32}, G_{33}, G_{34} \}$ with
\[
G_{31} = \left( \begin{array}{rr}
0   & \frac{1}{4}   \\[2mm]       1  & \frac{3}{4}
\end{array}
\right), \qquad
G_{32} = G_{31}, \qquad
G_{33} = \left( \begin{array}{rr}
1 & -\frac{1}{4}  \\[2mm]  0 & -\frac{1}{4}
\end{array}
\right), \qquad
G_{34} = \left( \begin{array}{rr}
-1 & \frac{1}{4} \\[2mm] -4 & \frac{3}{4}
\end{array}
\right).
\]
All previous families have joint spectral radius $1$.
The $L_1$-norm is extremal for $G_1$.
This means that $\|G_{1j}\|_1 \le 1$ for $j=1, \ldots , 4$.
Hence $\rho(G_1) = 1$.
Since $G_2 = G_1$, it follows that $\rho(G_2) = 1$.
For the family $G_3$, we apply Algorithm~\ref{algoP}
and obtain the invariant polytope~$P$. In this case $P$ is an octagon with vertices
\[
\left\{
\pm \left( \begin{array}{r} \frac{1}{4} \\ 1 \end{array} \right),
\pm \left( \begin{array}{r} 1 \\ 0 \end{array} \right),
\pm \left( \begin{array}{r} 0 \\ 1 \end{array} \right),
\pm \left( \begin{array}{r} 1 \\ 4 \end{array} \right).
\right\}
\]
Thus, $\rho(G_3) = 1$, and hence $\rho(\cA_2) = \max \{\rho(G_1), \rho(G_2), \rho(G_3)\} = 1$.

\smallskip

\textbf{The proof of Theorem~\ref{th40}}

\smallskip

We start with introducing some further notation.
A norm $\|\cdot \|$ in $\re^d$ is called {\em extremal} for a family~$\cA$
if $\|A_i\| \le \rho(\cA) $ for all $A_i \in \cA$. Algorithm~\ref{algoP}
constructs an extremal polytope norm. A family~$\cA$ is called {\em product bounded}
if norms of all products of matrices from $\cA$ are uniformly bounded 
(see e.g. \cite{GZ01}).
If a family has an extremal norm and $\rho(\cA) = 1$, then it is product bounded.

We have shown that $\rho(\cA_1)  = \rho(\cA_2) = 1$. Hence, the block lower-triangular form yields that $\rho(\cA) = \max \{\rho(\cA_1), \rho(\cA_2)\} = 1$, and so
$\rho \{T^{(3)}_i, \ i = 1, \ldots , 4\} = \frac14\, $.
Furthermore, all matrices $T_i^{(1)}$ in a special basis of the space~$\re^{24}$ have the form:
\begin{equation}\label{factor1}
T_i^{(1)} \ = \
\left(
\begin{array}{ccc}
J_2 & 0 & 0\\
* & J_3 & 0\\
* & * & T_i^{(3)}
\end{array}
\right)
\, , \qquad i = 1,2,3,4,
\end{equation}
where $J_s$ is the $(s+1)\times (s+1)$-diagonal matrix with all
diagonal entries equal to $2^{-s}$ (see~\cite{DL92, NPS}). Therefore,
the joint spectral radius of~$\{T^{(1)}_i, i = 1,\ldots , 4\}$
is equal to the maximum of the joint spectral radii of the three blocks, i.e.,
the maximum of $\rho(J_2) = \frac14$, of $\rho(J_3) = \frac18$, and
of~$\rho \{T^{(3)}_i, i = 1,\ldots , 4\} = \frac14$.
Thus, $\rho \{T^{(1)}_i, i = 1,\ldots , 4\} = \frac14$, and hence
$\alpha_{\varphi} = -\log_2 \frac14 = 2$. The H\"older exponent is found.
Now let us analyze the regularity of the derivative~$\varphi\,'$.

For any refinable function~$\varphi$, the  modulus of continuity 
$\omega_{\varphi\,'}(h)$ is asymptotically equivalent to the logarithm of the 
left-hand side of the equality~(\ref{defect1})
with $k = - [\log_2 h]$ (see~\cite{P06}). Hence, to prove that 
$\omega_{\varphi\,'}(h) \asymp h \, |\log h|^2$
it suffices to establish~(\ref{defect1}). Applying factorization~(\ref{factor1}) and
the results  of Steps 1-3, we obtain
$$
T_i^{(1)} \ = \
\left(
\begin{array}{cccc}
J_2 & 0 & 0 & 0\\
* & J_3 & 0 & 0\\
* & * & \frac14 B_i & 0\\
* & * & * & \frac14 C_i
\end{array}
\right)
\, , \qquad i = 1,2,3,4.
$$
In this block lower-triangular form, we have three blocks
($J_2, \frac14 B_i$ and $\frac14 C_i$) with the joint spectral radius $\frac14$ and
one ($J_3$) with a smaller spectral radius ($\frac18$). Moreover, all these former three blocks
are product bounded, since they have extremal norms. Therefore~\cite{P06},
\begin{equation}\label{defect2}
\max \, \|T_{i_k}^{(1)}\cdots T_{i_1}^{(1)}\| \ \le \ C_1 4^{-k}k^2 \, , \quad k \in \n\, ,
\end{equation}
where $C_1$ is a constant.
On the other hand, it is verified directly that each of the  matrices $T_i^{(1)}, i = 1,2,3$, has two  Jordan blocks
of size $3$ corresponding to the leading eigenvalue $\lambda = \frac14$. Hence,
the left-hand side of~(\ref{defect2}) is bigger than or equal to $\|[T_1^{(1)}]^k\| \ge C_2 \lambda^k k^2 = 4^{-k}k^2$.
Therefore, it is asymptotically equivalent to $4^{-k}k^2$. This proves~(\ref{defect1}) and hence
$\omega_{\varphi\,'}(h) \, \asymp \, h \, |\log h|^2$.
{\hfill $\Box$}
\smallskip

\subsection{Other values of the parameter $\omega$}

The convergence analysis of the Butterfly scheme can be
extended to other values of
$\omega \in [0,\frac{1}{4}]$. In this case we have to deal with
$24 \times 24$-matrices $T_i^{(1)}, i = 1,2,3,4.$
The scheme converges (to continuous limit functions) if and only if
their joint spectral radius $\rho$ is smaller than one.
The regularity of the scheme is equal to
$\alpha_{\varphi} = -\log_2 \rho$.

For  $\omega=\frac{1}{4}$, each matrix $T_i^{(1)}, i=1,2,3$, has two simple eigenvalues of modulus one: precisely $1$ and $-1$.
The  matrix $T_4^{(1)}$ has a simple eigenvalue $-1$ and the  $1$ of multiplicity~$2$
(both algebraic and geometric).
The two leading eigenvectors corresponding to $1$ and $-1$ define
a common invariant subspace for the family which can be transformed
into a similar block triangular form.
The $2 \times 2$-blocks are respectively
\[
\left( \begin{array}{rr}
1   & 0   \\[2mm]       0  & -1
\end{array}
\right), \qquad
\left( \begin{array}{rr}
-\frac12 & \frac{\sqrt3}2 \\[2mm]  \frac{\sqrt3}2 & \frac12
\end{array}
\right), \qquad
\left( \begin{array}{rr}
-\frac12 & -\frac{\sqrt3}2 \\[2mm]  -\frac{\sqrt3}2 & \frac12
\end{array}
\right), \qquad
\left( \begin{array}{rr}
1 & 0 \\[2mm] 0 & 1
\end{array}
\right).
\]
They are all symmetric, hence their joint spectral radius equals to the maximal spectral radius of these matrices~\cite{CH}, i.e., is equal to one.
The remaining $22 \times 22$ family of $4$ matrices has the fourth
matrix as an s.m.p. Starting from its (unique) leading eigenvector Algorithm~\ref{algoP} terminates within $8$ iterations and constructs an invariant polytope norm with $487$ vertices.
This proves that $\rho (T^{(1)}_1, \ldots , T^{(1)}_4) = 1$, and hence the scheme does not converge.
%
\begin{table}[h]
\begin{center}
\hskip -5mm
\begin{tabular}{|c|ccccccccc|}
\hline
$k$
    & $0$ & $1$ & $2$ & $3$ & $4$ & $5$ & $6$ & $7$ & $8$ \\
\hline
\rule{0pt}{9pt}\noindent
$\alpha_\varphi(\omega_k)$
    & $1$ & $1.1000$ & $1.2284$ & $1.4150$ & $2$ & $1.6781$ & $1.4150$ & $1.1926$ & $1$ \\
\hline
\hline
$k$
    & $9$ & $10$ & $11$ & $12$ & $13$ & $14$ & $15$ & $16$ & $ $ \\
\hline
\rule{0pt}{9pt}\noindent
$\alpha_\varphi(\omega_k)$
    & $0.8301$ & $0.6781$ & $0.5406$ & $0.4150$ & $0.2996$ & $0.1926$ & $0.0931$ & $0^-$ & $ $ \\
\hline
\end{tabular}\\[0.2mm]
\caption{Computed H\"older exponent of Butterfly scheme for
$\omega_k = \frac{k}{64}$, $k=0,\ldots,16$.\label{tab:omega}}
\end{center}
\end{table}

\begin{figure}[ht]
\centering \global\def\path{#1}\input{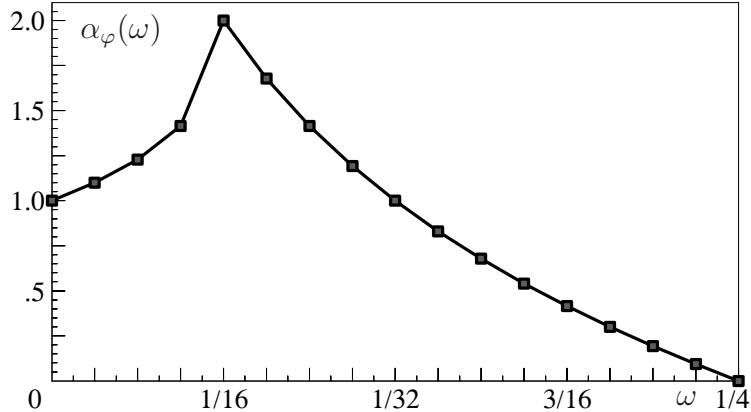}
\vspace{4mm}
\caption{The computed H\"older exponent of the generalized Butterfly scheme.
\label{fig:HolderGBS}}
\end{figure}

We have successfully applied our procedure also for other values
$\omega \in \left[ 0, \frac{1}{4} \right)$. This leads
us to conjecture that the generalized Butterfly subdivision scheme
is convergent in the whole interval. To support this conjecture we
report in Table \ref{tab:omega} the results obtained for $\omega = \frac{k}{64}$,
$k=0,1,\ldots,15,16$ (see also Figure \ref{fig:HolderGBS}).

\section{Applications: the regularity of Daubechies wavelets}\label{s6}

One of the most important applications of the joint spectral radius is the computation of the H\"older regularity of refinable functions and wavelets.
For Daubechies wavelets,  this problem was studied in many works 
(see~\cite{CD, CH, D88, DL92, G, NPS, Rio, Vil} and references therein).
Let us recall that the Daubechies wavelets is a system of functions $2^{j/2}
\psi (2^j x - n), \, j, n \in \z$, that constitutes an orthonormal basis in
$L_2(\re)$. All functions of this system are generated by
double dilates and integer translates of the compactly supported {\em wavelet function} $\psi$.
I.Daubechies in~\cite{D88} constructed a countable
family of wavelet functions~$\psi = \psi_N, \, N \ge 1$, each generates its own
 wavelet system. The function $\psi_1$ is the Haar function.
For all $N\ge 2$ the functions $\psi_N$ are continuous, their smoothness increases in $N$
and $\alpha_{\psi_N} > 0.2 N$~\cite{D88}. So, there are arbitrarily smooth systems of wavelets.
However, the price for the regularity is the length of the support, which also grows with~$N$:
${\rm supp}\, \psi_N\, = \, [0, 2N-1]$. The regularity  is a very important characteristics
of wavelets, in particular, for their applications in functional analysis, approximation theory, image processing and in numerical PDE. There are several methods to obtain lower and upper bounds for the H\"older exponents of wavelet functions (see~\cite{CD, D88, NPS, Rio, Vil}).
The matrix approach is the only one that theoretically allows to find them precisely.
It was established in~\cite{CH, DL92} that  $\alpha_{\psi_N} = N-\log_2 \rho(B_0, B_1)$,
where $B_0, B_1$ are special  matrices of size $(N-1)\times (N-1)$. This enabled
to find the precise values of the H\"older exponent for some small $N$.
For $N=2, 3$, and $4$, the value $\alpha_{\psi_N}$
were found by Daubechies and Lagarias in~\cite{DL92}; for $N=5, 6, 7$, and $8$, they were computed by G.Gripenberg~\cite{G}.
Every time a delicate analysis of special properties of those matrices was involved. In all the cases the s.m.p. of the family $\{B_0, B_1\}$ was one of those two matrices, and
it was a general belief that this is the case for all~$N$. We apply the standard
   routine of Algorithm~\ref{algoP} to find the precise values of $\alpha_{\psi_N}$ for all $N \le 20$.
   In particular, we shall see that for $N=10$, the conjecture of one matrix s.m.p. is violated and
   the s.m.p. is $B_0^2B_1^2$.

We need to  recall key steps of construction of the matrices~$B_0, B_1$.
For every $N=1,2,\ldots$ we have a set of $2 N$ {\em Daubechies filter coefficients}:
 $c_0, \ldots, c_{2N-1}$. They possess some special properties,  in particular,  
$\sum_{i=0}^{2N-1}c_i = 2$ and the polynomial  $m(z) = \sum_{n=0}^{2N-1} c_n z^n$ has 
zero of order $N$ at the point $z=-1$.
We set
$$
q(z)  =  \frac{m(z)}{\left((1+z)/2\right)^N} = \sum_{n=0}^{N-1} q_n z^n.
$$
and  write the transition $k \times k$-matrices
as follows:
\begin{eqnarray*}
&& (B_0)_{ij} =  q_{2i-j-1}, \qquad (B_1)_{ij} = q_{2i-j}, 
\qquad \qquad \qquad i, j  = 1,\ldots,N-1.
\end{eqnarray*}
We compute $\rho(B_0, B_1)$ by Algorithm \ref{algoP}. For some $N$ we have a non-unique s.m.p.
(these are the cases when $B_0$ and $B_1$ are both s.m.p.) and find the balancing vector
$\pa$ by the method in Subsection~\ref{s3.4}. However, due to symmetry, the entries of that
two-dimensional vector~$\pa$ are equal.  Another difficulty, much more significant for the
Daubechies matrices is that with growing~$N$ they become very ill-conditioned.
All vertices of the constructed polytope have very small last components, which
corresponds to the property of
quasi-invariance of a certain  subspace and determines a polytope strongly flattened
along certain directions.
For $N=10$, the last components  are about $10^{-12} - 10^{-19}$ of the values of the first components.
This creates enormous numerical difficulties in the running of Algorithm~\ref{algoP}, in particular, in the linear programming routines. That is why we use the technique with extra initial vertices 
(Section~\ref{s4}).
In the next subsections  we present three illustrative cases ($N=4,10,12$) and report the
computed H\"older regularity of Daubechies wavelets for al $N \le 20$.

\subsection{Illustrative examples}\label{s6.1}

We demonstrate the computation process for $N=4, 10$ and $12$ and see the crucial changes in the behaviour of Algorithm \ref{algoP} when the dimension grows. The case $N=4$ was done (by a different approach) by Daubechies and Lagarias in~\cite{DL92}, while the two other cases are new.

\subsection*{The case $N=4$.}

For the pair of transition matrices:
\begin{eqnarray*}
B_0 & = &
\left(
\begin{array}{rrr}
   5.212854848820774 &                  0 &                  0 \\
   1.703224934278843 & -4.676287953813834 &  5.212854848820774 \\
                   0 & -0.239791829285782 &  1.703224934278843
\end{array}
\right),
\\[3mm]
B_1 & = &
\left(
\begin{array}{rrr}
  -4.676287953813834 &  5.212854848820774 &                  0 \\
  -0.239791829285782 &  1.703224934278843 & -4.676287953813834 \\
                   0 &                  0 & -0.239791829285782
\end{array}
\right),
\end{eqnarray*}
the candidate s.m.p. is $B_0$ with $\rho(B_0) = 5.212854848820774\ldots$.
In order to apply Algorithm \ref{algoP} we compute the leading eigenvector
of $B_0$,
\[
v_1 = \left(
\begin{array}{r}
  1.0000      \\
  0.1662   \\
 -0.0113
\end{array}
\right)
\]
and set $\cV_0 = \{v_1 \}$.
Observe that $v_1$ almost lies on the subspace $E_2 \subset \re^3$ spanned by the
vectors $\{ e_1, e_2 \}$ of the canonical basis of~$\re^3$.
Applying the normalized matrices $\widetilde B_0$ and $\widetilde B_1$
repeatedly to $\cV_0$
one observes that the resulting vectors also almost lie on the
subspace $E_2$. This has implications on the flatness of the invariant
polyhedron computed by Algorithm \ref{algoP}.
\begin{figure}[ht]
\centering
\includegraphics[width=8.4cm]{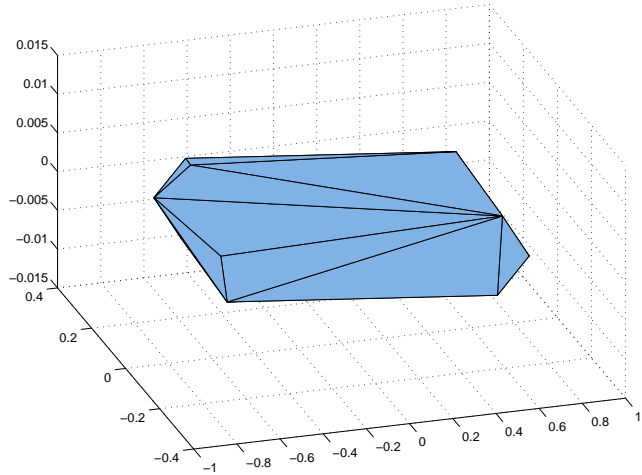}
\caption{The invariant polytope for the Daubechies matrices for $N=4$ computed
by the standard algorithm. \label{fig:D41}}
\end{figure}
In its basic implementation the algorithm converges and
generates a centrally symmetric polytope in $\re^3$ of $6 \cdot 2$
vertices, in $7$ iterations.
The partial polytope norms of the family $\widetilde \cA = \{ \widetilde B_0, \widetilde B_1 \}$ are reported in Table \ref{tab:k10}.
\begin{table}[h]
\begin{center}
\hskip -5mm
\begin{tabular}{|c|cccc|}\hline
Iteration $k$ & $4$ & $5$ & $6$ & $7$ \\
\hline
\rule{0pt}{9pt}\noindent
$\| \cdot \|_{\cP_i}$ &
$1.7845$ &  $1.1288$  &  $1.0571$ &  $1$
\\
\hline
\end{tabular}\\[0.1mm]
\caption{The partial polytope norms computed by the standard algorithm for $k=4$.
\label{tab:k4}}
\end{center}
\end{table}
The vertices (beyond $v_1$) of the polytope follow (we report only half of them):
\begin{eqnarray*}
&& v_2 = \left(
\begin{array}{r}
 -0.7308 \\
  0.0185 \\
  5.2249 \cdot 10^{-4}
\end{array}
\right), \quad
v_3 = \left(
\begin{array}{r}
 0.7308  \\
 0.2548  \\
 6.8062 \cdot 10^{-4}
\end{array}
\right), \quad
v_4 = \left(
\begin{array}{r}
 -0.7308 \\
 -0.0108 \\
  0.0115
\end{array}
\right),
\\[3mm]
&& v_5 = \left(
\begin{array}{r}
 -0.7308  \\
 -0.2175  \\
  0.0042
\end{array}
\right), \quad
v_6 = \left(
\begin{array}{r}
 -0.7308  \\
 -0.0393 \\
  0.0113
\end{array}
\right)
\end{eqnarray*}
and the corresponding unit polytope
$P \, = \, {\rm absco} \{  v_i, \  i  = 1, \ldots , 6\}$
is shown in Figure \ref{fig:D41}. We see that $P$ appears to be very flat.
%
\begin{figure}[hbt]
\centering
\includegraphics[width=8.4cm]{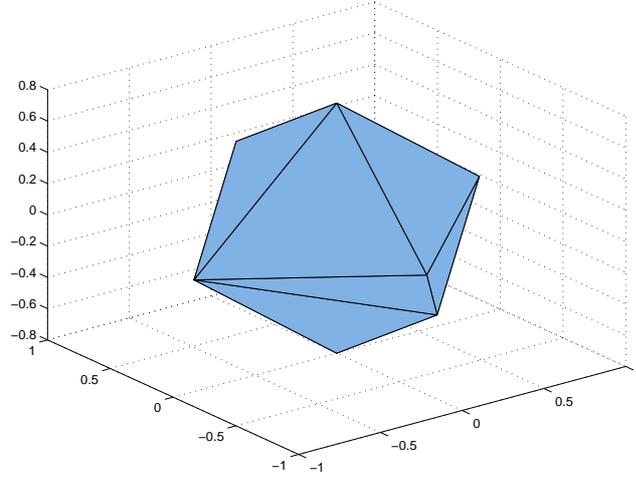}
\caption{Polytope extremal norm for the Daubechies matrices for $N=4$ computed
with an extra initial vertex. \label{fig:D42}}
\end{figure}

In fact, the largest singular value\footnote{Recall that for a matrix $B \in \R^{p,q}$
(or $B \in \co^{p,q}$), the {\em reduced} singular value decomposition is given
by $B = U \Sigma W^*$ where $U \in \co^{p,q}$ and $W \in \co^{q,q}$ are unitary
matrices and $\Sigma \in \R^{q,q}$ is a diagonal matrix with nonnegative diagonal elements 
$\{ \sigma_i \}_{i=1}^{q}$ (the singular values, usually ordered
in decreasing way) that are the square roots of the eigenvalues of the Hermitian
(semi)-positive definite matrix $B^* B$ (see e.g. \cite{GVL13}).}
of the matrix of its vertices $V = \{  \pm v_i, \  i  = 1, \ldots , 6\}$
is $\sigma_1 = 1.9365 $, while the smallest one is $\sigma_3 = 1.0923 \cdot 10^{-2}$
(note that if $\sigma_3$ would be zero then $V$ would not span the whole space and
the polytope would be contained in a subspace).
This almost $200$  times difference gives a numerical evidence of the flattening phenomenon.
To avoid it we add an extra initial vertex $x_1 = \frac45 e_3$
and obtain a better behaviour of Algorithm \ref{algoP} and a more balanced
polytope (see Figure \ref{fig:D42}).
The polytope has now $4 \cdot 2$ vertices which are computed in only two iterations.
Denote $v_2 = x_1$.
The vertices beyond $v_1$ and $v_2$ of the polytope follow (we report only half of them):
\[
v'_3 = \left(
\begin{array}{r}
            0 \\
  0.8000 \\
  0.2613
\end{array}	
\right), \qquad
v'_4 = \left(
\begin{array}{r}
            0 \\
 -0.7176 \\
 -0.0368
\end{array}
\right).
\]
The largest singular value is now  $\sigma_1 = 1.1597 \cdot 10^0$ and
the smallest singular value is now $\sigma_3 = 7.7402 \cdot 10^{-1}$, which demonstrate
a much more balanced shape of the unit ball of the polytope extremal norm.
The computed H\"older exponent of $\psi_{4}$ is
$\alpha_{\varphi_{4}} \ = \ 4 -\log_2 \rho(B_0, B_1) =
4 \, - \,  \log_2 \rho(B_0) \, = \, 1.6179\ldots$.

\subsection*{The case $N=10$.}

We have the $9\times 9$-matrices $B_0, B_1$ and are going to prove that
the s.m.p. is
\[
\Pi \ = \ B_0^2\,B_1^2, \qquad \rho_c = \rho(\Pi)^{1/4} = 99.636965469277555\ldots
\]
which  interesting in itself since it contradicts to the conjectured property
that $\rho( B_0, B_1 ) = \max\{ \rho(B_0), \rho(B_1) \}$ for all~$N$ (see Introduction).

Let $\widetilde B_0 = B_0/\rho_c$, $\widetilde B_1 = B_1/\rho_c$.
Applying Algorithm \ref{algoP} we compute the starting set of vectors
$\cV_0 = \{ v_1, v_2, v_3, v_4 \}$, where
$v_1$ is the leading eigenvector of $\Pi$, $v_2 = \tilde B_1 v_1, \, v_3 = \tilde B_1 v_2, v_4 = \tilde B_0 v_3$.
Clearly, $\tilde B_0v_4 = v_1$. Thus, $v_2, v_3$ and $v_4$ are the leading eigenvectors of cyclic permutations of the product~$\Pi$. We have
\begin{eqnarray*}
& v_1 =
\left(
\begin{array}{r}
  1.7122 \cdot 10^{-3} \\
  1.0000 \cdot 10^{0\ } \\
  3.0340 \cdot 10^{-1} \\
 -3.0515 \cdot 10^{-1} \\
 -5.9219 \cdot 10^{-2} \\
  5.9518 \cdot 10^{-4} \\
  9.0971 \cdot 10^{-5} \\
 -2.5228 \cdot 10^{-7} \\
 -6.6280 \cdot 10^{-12}
\end{array}
\right), \qquad
& v_2 =
\left(
\begin{array}{r}
  3.8568 \cdot 10^{-1} \\
  1.1649 \cdot 10^{0\ } \\
  3.6053 \cdot 10^{-2} \\
 -1.8840 \cdot 10^{-1} \\
 -1.8276 \cdot 10^{-2} \\
  5.0899 \cdot 10^{-4} \\
  1.0986 \cdot 10^{-5} \\
 -1.8367 \cdot 10^{-8} \\
  1.2778 \cdot 10^{-15}
\end{array}
\right),
\end{eqnarray*}
\begin{eqnarray*}
& v_3 =
\left(
\begin{array}{r}
  1.1395 \cdot 10^{-2} \\
  1.2195 \cdot 10^{0\ } \\
  5.6969 \cdot 10^{-1} \\
 -2.4167 \cdot 10^{-2} \\
 -1.2787 \cdot 10^{-2} \\
 -9.4332 \cdot 10^{-5} \\
  4.9101 \cdot 10^{-6} \\
 -2.1784 \cdot 10^{-9} \\
 -2.4634 \cdot 10^{-19}
\end{array}
\right), \qquad
& v_4 =
\left(
\begin{array}{r}
  4.4172 \cdot 10^{-3} \\
 -1.1524 \cdot 10^{0\ }  \\
 -8.3168 \cdot 10^{-1} \\
  4.6863 \cdot 10^{-2} \\
  3.7374 \cdot 10^{-2} \\
  7.1841 \cdot 10^{-4} \\
 -3.4331 \cdot 10^{-5} \\
  3.4387 \cdot 10^{-8} \\
  4.1997 \cdot 10^{-13}
\end{array}
\right).
\end{eqnarray*}
Observe that the last two components of all the vectors~$v_1, \ldots , v_4$ are very small, 
i.e., all these vectors almost lie on the subspace $E_7 \subset \re^9$ spanned by the
vectors $\{ e_1, \ldots, e_7 \}$ of the canonical basis of~$\re^9$.
Applying $\widetilde B_0$ and $\widetilde B_1$ repeatedly to $\cV_0$
one observes that the resulting vectors also almost lie on the
subspace $E_7$. This means that the invariant polytope computed by
Algorithm \ref{algoP} is nearly degenerate, it is close to a $7$-dimensional polytope.
A consequence of this is a slow convergence behaviour of the algorithm and an 
ill-conditioning of basic linear algebra operations.
The algorithm terminates and
generates a centrally symmetric polytope in $\re^9$ of $220 \cdot 2$
vertices, in $16$ iterations.
The partial polytope norms of the family $\widetilde \cA = \{ \widetilde B_0, \widetilde B_1 \}$
are reported in Table \ref{tab:k10}.
\begin{table}[h]
\begin{center}
\hskip -5mm
\begin{tabular}{|c|cccccccccccc|}\hline
$k$ & $5$ & $6$ & $7$ & $8$ & $9$ & $10$ & $11$ & $12$ & $13$ & $14$ & $15$ & $16$ \\
\hline
\rule{0pt}{9pt}\noindent
$\| \cdot \|_{P_i}$ &
$24.856$ &  $4.693$  &  $2.990$ & $2.237$  &  $1.743$  &  $1.414$  &
$1.140$  &  $1.064$  &  $1.027$ & $1.025$  &  $1.001$  &  1
\\
\hline
\end{tabular}\\[0.1mm]
\caption{The partial polytope norms computed by the Algorithm \ref{algoP}
for $N=10$.\label{tab:k10}}
\end{center}
\end{table}
The largest singular value of the set of vertices is $\sigma_1 = 1.3253 \cdot 10^1$ and
the smallest singular value is $\sigma_9 = 1.0620 \cdot 10^{-10}$, which gives a numerical evidence of the flatness of the polytope. To correct the behavior of the algorithm we add an extra vector
along the ``most narrow'' (for the polytope~$P$) direction (see Section~\ref{s4}).
Adding the vector $v_5 = \frac12 e_9$, we obtain a better behaviour of the algorithm and a more balanced polytope.
The results are summarized here; the polytope has $75 \cdot 2$ vertices which are computed in $10$ iterations.
\begin{table}[h]
\begin{center}
\hskip -5mm
\begin{tabular}{|c|ccccccccc|}\hline
$i$ & $2$ & $3$ & $4$ & $5$ & $6$ & $7$ & $8$ & $9$ & $10$ \\
\hline
\rule{0pt}{9pt}\noindent
$\| \cdot \|_{\cP_i}$ &
$9.4019$  &  $1.8655$  &  $1.3861$  & $1.2622$  &
$1.2076$  &  $1.1150$  &  $1.0512$  &  $1.0169$ &  1
\\
\hline
\end{tabular}\\[0.1mm]
\caption{The partial polytope norms computed by the modified algorithm
for $N=10$.\label{tab:k10m}}
\end{center}
\end{table}
The largest singular value is now  $\sigma_1 = 7.7943 \cdot 10^0$ and
the smallest one is $\sigma_9 = 1.1401 \cdot 10^{-3}$, which demonstrate
a much more balanced shape of the unit ball of the polytope extremal norm.
The computed H\"older exponent of $\, \psi_{10}\, $ is
$$
\alpha_{\varphi_{10}} \ = \ 10\, - \, \log_2 \, \rho(B_0, B_1) \ =\
10 \, - \, \frac14 \, \log_2 \rho(B_0^2B_1^2) \ = \ 3.361390821401114\ldots
$$
\begin{remark}\label{r20b}
{\em If we consider the adjoint family $\{ B_0^{*}, B_1^* \}$
we have naturally $\rho(B_0, B_1) = \rho(B_0^*, B_1^*)$.
It is remarkable that applying Algorithm~\ref{algoP} to $\{ B_0^{*}, B_1^* \}$,
we do not have problems with flattening. Alas convergence remains slow
($12$ iterations and $370 \cdot 2$ vertices). }
\end{remark}

\subsection*{The case $N=12$.}

In this case  both $B_0$ and $B_1$ are s.m.p.,
i.e. $\rho( B_0, B_1 ) =  \rho(B_0) =  \rho(B_1)$.
The balancing technique (Section~\ref{s3}) gives $\alpha = (1,1)$, i.e.,  equal weights to the leading
eigenvectors of $B_0$ and of $B_1$. The two vectors, say $v_1$ and $v_2$ follow:
\begin{eqnarray*}
& v_1 =
\left(
\begin{array}{r}
            0 \\
  1.3465 \cdot 10^{-1} \\
  1.0000 \cdot 10^{0\ } \\
  4.3937 \cdot 10^{-1} \\
 -1.1888 \cdot 10^{-1} \\
 -4.1549 \cdot 10^{-2} \\
 -5.5942 \cdot 10^{-4} \\
  1.2299 \cdot 10^{-4} \\
  2.0148 \cdot 10^{-7} \\
 -3.6435 \cdot 10^{-9} \\
  1.1240 \cdot 10^{-13}
\end{array}
\right), \qquad
& v_2 =
\left(
\begin{array}{r}
 -4.8598 \cdot 10^{-1} \\
 -3.2838 \cdot 10^{-2} \\
  1.0000 \cdot 10^{0\ } \\
  3.5198 \cdot 10^{-1} \\
 -8.3828 \cdot 10^{-3} \\
 -5.6206 \cdot 10^{-3} \\
 -5.3176 \cdot 10^{-5} \\
  3.0075 \cdot 10^{-6} \\
 -3.7539 \cdot 10^{-9} \\
 -1.2233 \cdot 10^{-13} \\
            0
\end{array}
\right).
\end{eqnarray*}
We observe that they almost lie on the subspace $E_8 \subset \re^{11}$ spanned by the
vectors $\{ e_1, \ldots, e_{8} \}$  of the canonical basis of~$\re^{11}$.
Applying $\widetilde B_0$ and $\widetilde B_1$ repeatedly to $\cV_0$
one observes that the resulting vectors also almost lie on the
subspace $E_8$. Note that the last components seems to vanish
exponentially in the number of iterations (from the smallest to the highest index).
If we add an extra initial vector $v_3 = \frac12 e_{11}$, Algorithm \ref{algoP}
terminates after $9$ iterations with an invariant polytope  of $48 \cdot 2$ vertices.
The partial polytope norms are reported in Table \ref{tab:k12}.
The largest singular value  is $\sigma_1 = 3.8190 \cdot 10^0$ and
the smallest one is $\sigma_9 = 1.5416 \cdot 10^{-5}$, which demonstrate
a relatively balanced shape. The computed H\"older exponent is 
$\alpha_{\psi_{12}} = 12 - \log_2 \rho(B_0) = 3.833483495658518\ldots$.
\begin{table}[h]
\begin{center}
\hskip -5mm
\begin{tabular}{|c|cccc|}\hline
$i$ & $4$ & $5$ & $6$ & $7$ \\
\hline
\rule{0pt}{9pt}\noindent
$\| \cdot \|_{\cP_i}$ &
$5.0579$ &  $1.5497$  &  $1.1597$  & $1$
\\
\hline
\end{tabular}
\caption{The partial polytope norms computed by the modified algorithm
for $N=12$.\label{tab:k12}}
\end{center}
\end{table}
\vspace{-9mm}
\subsection{The table of results for $N \le 20$}\label{s6.2}

Proceeding this way we have computed the exact values of H\"older exponents of Daubechies wavelets according to Table \ref{tab:results}.

We indicate the s.m.p., the extra initial vectors,
the number of vertices of the final polytope ($\#V$), the number of iterations of Algorithm 
\ref{algoP} ($\#$its) and the H\"older exponent $\alpha$.
\begin{table}[h]
\begin{center}
\hskip -5mm
\begin{tabular}{|c|ccccc|}\hline
$N$ & s.m.p.  &  Extra vertices  &  $\#$its  & $\#V$  &  $\alpha$ \\
\hline
\rule{0pt}{9pt}\noindent
$2$  & $B_0$ & {\rm none} & $1$ & $1 \cdot 2$ & $0.55001\ldots$ \\
$3$  & $B_0$ & {\rm none} & $3$ & $3 \cdot 2$ & $1.08783\ldots$ \\
$4$  & $B_0$ & $ 0.8 e_{3} $ & $2$ & $4 \cdot 2$ & $1.61792\ldots$ \\
$5$  & $B_0\ \mbox{and} \ B_1$ & $ 0.1 e_{4} $ & $4$ & $8 \cdot 2$ & $1.96896\ldots$ \\
$6$  & $B_0\ \mbox{and} \ B_1$ & $ 0.1 e_{5} $ & $5$ & $11 \cdot 2$ & $2.18913\ldots$ \\
$7$  & $B_0\ \mbox{and} \ B_1$ & $ 0.1 e_{5} $ & $5$ & $12 \cdot 2$ &
$2.46040\ldots$ \\
$8$  & $B_0\ \mbox{and} \ B_1$ & $ 0.1 e_{7} $ & $5$ & $18 \cdot 2$ &
$2.76081 \ldots$ \\
$9$  & $B_0\ \mbox{and} \ B_1$ & $ 0.5 e_{8} $ & $6$ & $24 \cdot 2$ &
$3.07361\ldots$ \\
$10$ & $B_0^2 B_1^2$ & $ 0.5 e_{9} $ & $10$ & $90 \cdot 2$ & $3.36139\ldots$ \\
$11$ & $B_0\ \mbox{and} \ B_1$ & $0.5 e_{10}$ & $11$ & $75 \cdot 2$ & $3.60346\ldots$ \\
$12$ & $B_0\ \mbox{and} \ B_1$ & $0.5 e_{11}$ & $7$ & $48 \cdot 2$ & $3.83348\ldots$ \\
$13$ & $B_0\ \mbox{and} \ B_1$ & $e_{12}$ & $18$ & $73 \cdot 2$ & $4.07347\ldots$ \\
$14$ & $B_0\ \mbox{and} \ B_1$ & $0.5 e_{13}, 0.25 e_{12} $ & $15$ & $73 \cdot 2$ &
$4.31676\ldots$ \\
$15$ & $B_0^4 B_1^2$ & $10^{-3} \{ e_{k} \}_{k=9}^{14}$ & $14$ & $376 \cdot 2$ &
$4.55611\ldots$ \\
$16$ & $B_0^2 B_1^2$ & $10^{-2} \{ e_{k} \}_{k=11}^{15} $ & $13$ & $372 \cdot 2$ & $4.78643\ldots$ \\
$17$ & $B_0\ \mbox{and} \ B_1$ & $10^{-3} \{ e_{k} \}_{k=11}^{16}$ & $11$ & $480 \cdot 2$ & $5.02444\ldots$ \\
$18$ & $B_0\ \mbox{and} \  B_1$ & $10^{-3} \{ e_{k} \}_{k=12}^{17}$ & $13$ & $409 \cdot 2$ &
$5.23915\ldots$ \\
$19$ & $B_0\ \mbox{and} \ B_1$ & $10^{-3} \{ e_{k} \}_{k=13}^{18}$ & $19$ &
$1395 \cdot 2$ & $5.46529\ldots$ \\
$20$ & $B_0\ \mbox{and} \ B_1$ & $10^{-3} \{ e_{k} \}_{k=13}^{19}$ & $24$ &
$2480 \cdot 2$ & $5.69116\ldots$ \\
\hline
\end{tabular} 
\caption{Computed H\"older exponent of Daubechies wavelets.\label{tab:results}}
\end{center}
\end{table}
%
\medskip

\section*{Appendix 1. Proof of Theorem~\ref{th10}.}
%
Necessity. If the algorithm terminates within finite time, then the products $\{ \Pi_i\}_{i=1}^r$ are dominant and their leading eigenvalues are unique and simple. This is shown in the same way as
in the proof of Theorem~4 of~\cite{GP13} for $r=1$. To prove that $\pa$ is admissible, we take arbitrary $i$ and $j$ and denote by $z$ the vertex of the final polytope~$P$ with the largest scalar product $(v_j^*, z)$.
Since $(v_j^*, \tilde \Pi_j^kz) = ([\tilde \Pi_j^*]^kv_j^*,z) = (v_j^*, z)$, all the points $\{\tilde \Pi_j^kz\}_{k \in \n}$ also provide the largest scalar product with the vector $v_j^*$.
Hence, they are all on the boundary of $P$, i.e., they are not absorbed in the algorithm. Consequently, the algorithm can terminate within finite time only if
$z$ is the leading eigenvector of $\tilde \Pi_j$, i.e.,  $z = \alpha_jv_j$. Thus, the maximal scalar product $(v_j^*, z)$ over all $z \in P$ is attained at a unique vertex $z = \alpha_jv_j$, where it 
is equal to  $(v_j^*, \alpha_j v_j) = \alpha_j$. Since $i \ne j$, it follows that
$$\sup\limits_{z \in \alpha_i P_{i, \infty}}(v_j^*, z) < \alpha_j.$$ 
Thus, $\alpha_i q_{ij} < \alpha_j$, which proves the admissibility of $\pa$.

\smallskip

Sufficiency. Denote by $\Omega$ the set of products $\tilde \Pi_i, \, i = 1, \ldots , r$, and of its cyclic permutations.
Since this is a set of dominant products for $\tilde \cA$, their leading
eigenvectors~$\{v_i^{(k)} \ | \  k = 1, \ldots , n_i\, , \  i = 1, \ldots , r\}$
 are all different up to normalization, i.e., they are all non-collinear. Indeed, if,
say, $v_i^{(k)} = \lambda\, v_j^{(l)}\, , \ \lambda \ne 0$, then, replacing the products
$\Pi_i$ and $\Pi_j$ by the corresponding cyclic permutations, it may be assumed that
$k = l = 1$. However, in this case $\tilde \Pi_j^{k_2} \tilde \Pi_i^{k_1}v_i^{(1)} =
\tilde \Pi_j^{k_2}\, v_i^{(1)} = v_i^{(1)}$
for any $k_1, k_2$. Therefore,  the spectral radius of every product of the form $\tilde \Pi_j^{k_2} \tilde \Pi_i^{k_1}$
is at lest one. By the dominance assumption, this product is a power of some product
$\tilde \Pi \in \Omega$. Taking now $k_1, k_2$ large enough and applying Lemma~\ref{l30}
first to the words $a = \tilde \Pi, b = \tilde \Pi_j$ and then to
 the words $a = \tilde \Pi, b = \tilde \Pi_i$, we conclude that both $\tilde \Pi_j, \tilde \Pi_i$
 must be cyclic permutations of $\tilde \Pi$, which is impossible. Thus, all the leading 
eigenvectors~$\{v_i^{(k)}\}$ are non-collinear. Hence, there is $\varepsilon > 0$ such that for 
every $x \in \re^d \setminus \{0\}$ the ball of radius $\varepsilon \|x\|$ centered at~$x$ may contain leading eigenvectors of at most one matrix from~$\Omega$.

If the polytope algorithm with the initial roots~$\, \alpha_1  \cH_1, \ldots , \alpha_r \cH_r$ does not converge, then there is an element of some root, say, 
$\alpha_1v_1 = \alpha_1 v_1^{(1)} \in \alpha_1\cH_1$ and an infinite sequence 
$\{\tilde A_{b_k}\}_{k\in \n}$, which is not periodic with period $\tilde \Pi_1$,  
and such that every vector $u_k = \tilde A_{b_k-1}\cdots \tilde A_{b_1} \alpha_1 v_1$ 
is not absorbed in the algorithm.
This implies that there is a constant $C_0 > 0$ such that $\|u_k\| \ge C_0$ for all $k$.
On the other hand, $\|u_k\| \, \le \, M \, \|\alpha_1\, v_1\|$, hence
the compactness argument yields the existence of a limit point~$u \ne 0$
of this sequence. Thus, for some subsequence, we have $u_{j_k}\to u$ as $k \to \infty$.
Let $\delta > 0$ be a small number to specify. 
Passing to a subsequence, it may be assumed that $\|u_{j_n} -u_{j_k}\| \le \delta $
for all~$k, n$. Denote $G_k = \tilde A_{j_{k+1}-1}\cdots \tilde A_{j_k}$.
We have $u_{j_{k+1}} = G_k u_{j_k}\, , \ k \in \n$.
Invoking the triangle inequality, we obtain
$$
\bigl\| G_k u - u\bigr\|  \ \le \ \bigl\| G_k(u-u_{j_k})\bigr\|\, + \,
  \bigl\| G_k u_{j_k} - u_{j_k}\bigr\|\, + \,  \bigl\|u_{j_k} - u\bigr\|\ \le \
  M \delta + \delta + \delta \ = \ (M+2)\delta \, .
$$
Hence, Lemma~\ref{l20} yields $\rho(P_k) \ge 1 - C(d)M^{1+\frac{1}{d}} \delta^{\, 1/d}$ for all~$k$.
The dominance assumption implies that if $\delta > 0$ is small enough, then
all $G_k$ must be powers of matrices from  $\Omega$. Each matrix from $\Omega$ has a simple
unique leading eigenvalue~$1$. Therefore, there is a function $\mu (t)$
such that $\mu(t) \to 0 $ as $t \to 0$,
and for every matrix $Q$ which is a power of a matrix from $\Omega$
the inequality $\bigl\| Q u - u\bigr\| < t$ implies that there exists a leading eigenvector $w$ of $Q$
such that $\|w - u \| < \mu (t)$. Thus, for every $k \in \n$ there exists a leading eigenvector $w_k$ of $G_k$
such that $\|w_k - u\| < \mu \bigl( (M+2) \delta\bigr)$. For sufficiently small
$\delta$ we have $\mu \bigl( (M+2) \delta\bigr) < \varepsilon \|u\|$. Hence, for every $k$, the vector $w_k$ belongs to the ball of radius $\varepsilon \|u\|$ centered at~$u$. However, this ball
may contain a leading eigenvector of at most one matrix from~$\Omega$, say $\tilde \Pi$. Therefore, all $G_k, \, k \in \n$, are powers of $\tilde \Pi$ and $u$ is the leading eigenvector of $\tilde \Pi$. Thus, $u_{j_k} = \tilde \Pi^{\, p_k}\alpha_1v_1$ for some 
$p_k \in \n$.
Clearly,  $\tilde \Pi$ is a cyclic permutation of some $\tilde \Pi_j$. If $j \ne 1$, 
then assuming that $\tilde \Pi = \tilde \Pi_j$ (the general case is considered in the same way), we have $u_{j_k} \to \alpha_1 (v_j^{*}, v_1) v_j$ as $k \to \infty$. Since the balancing vector $\alpha$ is admissible and $(v_j^{*}, v_1) \le q_{1 j}$,
it follows that $\alpha_1 (v_j^{*}, v_1) < \alpha_j$. Therefore, the limit point $u = \alpha_1 (v_j^{*}, v_1) v_j = \lambda \alpha_jv_j$ for some $\lambda \in (0,1)$, is interior for the initial polytope ${\rm co } (\alpha \cH )$. 
This means that for large $k$, the point $u_{j_k}$ will be absorbed in the algorithm,
which contradicts to the assumption.
Consider the last case, when  $\tilde \Pi$ is a cyclic permutation of $\tilde \Pi_1$.
Since $u$ is the leading eigenvector of $\tilde \Pi$ we see that
$u = \beta v_1^{(s)}$ for some $s = 1, \ldots , n_1$ and $\beta \in \re$.
We assume $\beta > 0$, the case of negative $\beta$ is considered in the same way.
If $\beta  < 1$, then we again conclude that $u_{j_k}$ are absorbed in the algorithm
for  large $k$. If $\beta  \ge 1$, then for the product 
$\tilde \Pi_0 = \tilde A_{d_{n}}\cdots \tilde A_{d_{s}}$, we
have $\tilde \Pi_0 u = \beta v_1$, and hence $\tilde \Pi_0\tilde \Pi^k v_1 \to \beta v_1$.
Therefore, $\|\beta^{-1} \tilde \Pi_0\tilde \Pi^kv_1 \, - \, v_1\|\, \to \, 0$ as $k \to \infty$.
By Lemma~\ref{l20}, this means that the spectral radius of the product $\tilde \Pi_0\tilde \Pi^k$ tends to $\beta \ge 1$ as $k \to \infty$. The dominance assumption implies now that for every sufficiently large $k$, this product is a power of some $\tilde \Pi_a \in \Omega$. Applying Lemma~\ref{l30} to the words $a = \tilde \Pi_a, b = \tilde \Pi$, we see that $\Pi_a$ is a cyclic permutation of $\tilde \Pi$. 
In particular, $|\tilde \Pi_a| =  |\tilde \Pi| = n $. Therefore, the length 
$|\tilde \Pi_0\tilde \Pi^k| = (n-s+1) + kn = (k+1)n - (s-1)$ is divisible by 
$|\tilde \Pi_a| = n$. Hence, the number $(s-1)$ is divisible by $n$, which is impossible, because $s \le n$. 
This completes the proof.
{\hfill $\Box$}
\smallskip

\section*{Appendix 2: the Butterfly matrices.}

We report here the matrices of the Butterfly scheme for $\omega=\frac{1}{16}$.
First the $17 \times 17$ family $\cA$, then the $11 \times 11$ family $\cB$ and
finally the $6 \times 6$ family $\cC$.

\newpage

\begin{landscape}
The four matrices $A_1, A_2, A_3, A_4$:
\tiny
\begin{eqnarray*}
A_1 & = &
\left(
\begin{array}{rrrrrrrrrrrrrrrrr}
 \frac{54187}{4104} & \frac{22307}{4104} & \frac{131}{456} & -\frac{10201}{1368} & -\frac{15943}{4104} & \frac{8251}{4104} & \frac{3299}{1368} & -\frac{2399}{684} &
   \frac{8869}{1368} & \frac{27103}{4104} & \frac{18013}{4104} & -\frac{7727}{4104} & -\frac{1870}{513} & -\frac{4225}{2052} & -\frac{6887}{4104} & \frac{1099}{456}
   & \frac{887}{1368} \\[2mm]
 \frac{15049}{4104} & \frac{13997}{1026} & \frac{89}{342} & -\frac{4225}{1368} & -\frac{3062}{513} & \frac{7531}{4104} & \frac{120}{19} & -\frac{3851}{1368} &
   \frac{785}{456} & \frac{15553}{4104} & \frac{7957}{1026} & -\frac{9767}{4104} & -\frac{8171}{4104} & -\frac{449}{216} & -\frac{17}{54} & \frac{47}{171} &
   \frac{3911}{1368} \\[2mm]
 -\frac{29369}{4104} & -\frac{16895}{2052} & -\frac{773}{228} & \frac{4133}{1368} & \frac{475}{108} & \frac{5893}{4104} & -\frac{980}{171} & \frac{1643}{1368} &
   -\frac{4919}{1368} & -\frac{15629}{4104} & -\frac{7639}{2052} & \frac{5251}{4104} & \frac{9871}{4104} & \frac{15877}{4104} & \frac{7631}{2052} & -\frac{443}{228}
   & -\frac{2347}{1368} \\[2mm]
 \frac{65953}{4104} & -\frac{17155}{4104} & \frac{1475}{1368} & -\frac{9643}{1368} & \frac{1415}{4104} & \frac{133}{216} & -\frac{339}{152} & -\frac{1573}{684} &
   \frac{3469}{456} & \frac{30277}{4104} & \frac{1795}{4104} & \frac{1897}{4104} & -\frac{2441}{1026} & -\frac{299}{1026} & -\frac{7655}{4104} & \frac{5005}{1368} &
   -\frac{2347}{1368} \\[2mm]
 -\frac{43913}{4104} & \frac{34169}{4104} & -\frac{1283}{456} & \frac{4643}{1368} & -\frac{2677}{4104} & \frac{10537}{4104} & \frac{209}{72} & \frac{409}{684} &
   -\frac{6971}{1368} & -\frac{14153}{4104} & \frac{14935}{4104} & -\frac{7043}{4104} & \frac{1547}{2052} & \frac{4181}{2052} & \frac{14425}{4104} &
   -\frac{1429}{456} & \frac{3911}{1368} \\[2mm]
 -\frac{41957}{4104} & -\frac{1319}{2052} & -\frac{863}{171} & \frac{4577}{1368} & \frac{4501}{2052} & \frac{14623}{4104} & -\frac{42}{19} & \frac{1525}{1368} &
   -\frac{677}{152} & -\frac{19439}{4104} & -\frac{391}{513} & -\frac{5069}{4104} & \frac{4825}{4104} & \frac{12691}{4104} & \frac{4303}{1026} & -\frac{523}{171} &
   \frac{887}{1368} \\[2mm]
 \frac{889}{513} & -\frac{5081}{1026} & \frac{181}{1140} & \frac{112}{171} & \frac{3661}{1026} & \frac{991}{5130} & -\frac{1571}{855} & \frac{1309}{684} &
   \frac{1657}{1140} & -\frac{61}{540} & -\frac{7141}{2565} & \frac{1139}{2565} & -\frac{1963}{10260} & \frac{8447}{10260} & -\frac{977}{5130} & \frac{3871}{3420} &
   -\frac{191}{855} \\[2mm]
 -\frac{929}{270} & \frac{2515}{1026} & \frac{679}{380} & \frac{365}{342} & -\frac{3935}{1026} & -\frac{7168}{2565} & \frac{5017}{3420} & -\frac{1189}{855} &
   -\frac{599}{228} & -\frac{2255}{2052} & \frac{3863}{10260} & \frac{965}{1026} & \frac{13877}{10260} & -\frac{1199}{513} & -\frac{3193}{2565} & -\frac{3107}{3420}
   & -\frac{191}{855} \\[2mm]
 -\frac{1003}{10260} & -\frac{11509}{5130} & \frac{1411}{1140} & \frac{167}{180} & \frac{4391}{5130} & -\frac{13723}{10260} & \frac{1061}{1710} & \frac{3181}{1710}
   & \frac{401}{285} & -\frac{17719}{10260} & -\frac{7006}{2565} & -\frac{5029}{10260} & -\frac{8209}{10260} & -\frac{3826}{2565} & -\frac{20899}{10260} &
   \frac{53}{171} & -\frac{106}{171} \\[2mm]
 -\frac{8518}{2565} & -\frac{42173}{10260} & \frac{1021}{1140} & \frac{784}{285} & \frac{21967}{10260} & -\frac{533}{270} & -\frac{466}{171} & \frac{107}{228} &
   -\frac{10439}{3420} & -\frac{25}{1026} & -\frac{12499}{10260} & \frac{2981}{1026} & \frac{7787}{2565} & \frac{20039}{10260} & \frac{2999}{2565} &
   -\frac{115}{342} & -\frac{2429}{3420} \\[2mm]
 -\frac{56513}{5130} & -\frac{163009}{20520} & -\frac{11933}{6840} & \frac{5249}{1140} & \frac{66761}{20520} & -\frac{2971}{2565} & -\frac{1481}{360} &
   \frac{14827}{6840} & -\frac{2009}{380} & -\frac{14051}{2052} & -\frac{118553}{20520} & \frac{2269}{2052} & \frac{9701}{4104} & \frac{39097}{20520} &
   \frac{32347}{20520} & -\frac{2169}{760} & -\frac{1658}{855} \\[2mm]
 -\frac{51289}{10260} & \frac{51557}{5130} & -\frac{21}{95} & \frac{75}{76} & -\frac{12344}{2565} & \frac{3443}{10260} & \frac{18733}{3420} & -\frac{29}{57} &
   -\frac{1189}{684} & -\frac{20527}{10260} & \frac{6691}{2052} & -\frac{26377}{10260} & -\frac{9703}{10260} & -\frac{13943}{5130} & -\frac{793}{1026} &
   -\frac{5467}{3420} & \frac{775}{342} \\[2mm]
 -\frac{23953}{10260} & \frac{33437}{5130} & -\frac{523}{855} & \frac{509}{1140} & -\frac{3284}{2565} & \frac{2863}{2052} & \frac{688}{171} & \frac{3233}{3420} &
   \frac{47}{380} & -\frac{14779}{10260} & \frac{9307}{5130} & -\frac{6526}{2565} & -\frac{4882}{2565} & -\frac{2609}{2052} & -\frac{2063}{5130} & -\frac{761}{1140}
   & \frac{775}{342} \\[2mm]
 \frac{92413}{5130} & -\frac{132199}{20520} & \frac{3991}{2280} & -\frac{7453}{1140} & \frac{35951}{20520} & \frac{205}{513} & -\frac{11701}{6840} &
   -\frac{179}{760} & \frac{6845}{684} & \frac{65351}{10260} & -\frac{40487}{20520} & -\frac{5029}{10260} & -\frac{15133}{4104} & -\frac{2051}{1080} &
   -\frac{4313}{1080} & \frac{6557}{1368} & -\frac{1658}{855} \\[2mm]
 \frac{11347}{5130} & -\frac{3617}{2565} & \frac{4129}{1710} & -\frac{43}{285} & -\frac{151}{270} & -\frac{5068}{2565} & -\frac{3253}{3420} & -\frac{2231}{1710} &
   -\frac{109}{380} & \frac{22733}{10260} & \frac{1093}{2052} & \frac{13849}{5130} & \frac{19073}{10260} & \frac{415}{2052} & \frac{160}{513} & \frac{1193}{1140} &
   -\frac{2429}{3420} \\[2mm]
 -\frac{30997}{4104} & -\frac{79909}{20520} & \frac{2113}{1368} & \frac{42259}{6840} & \frac{78761}{20520} & -\frac{8299}{4104} & -\frac{2611}{2280} &
   \frac{1859}{570} & -\frac{29269}{6840} & -\frac{77327}{20520} & -\frac{52403}{20520} & \frac{48463}{20520} & \frac{32657}{10260} & \frac{8377}{5130} &
   \frac{19687}{20520} & -\frac{69}{760} & -\frac{637}{6840} \\[2mm]
 \frac{44899}{20520} & -\frac{21365}{4104} & \frac{9313}{6840} & -\frac{193}{1368} & \frac{403}{216} & -\frac{7219}{4104} & -\frac{4751}{6840} & \frac{6211}{3420} &
   \frac{13183}{6840} & -\frac{22907}{20520} & -\frac{84623}{20520} & -\frac{4067}{20520} & -\frac{12139}{10260} & -\frac{7481}{5130} & -\frac{54053}{20520} &
   \frac{7027}{6840} & -\frac{7129}{6840}
\end{array}
\right)
\\[2mm]
A_2 & = &
\left(
\begin{array}{rrrrrrrrrrrrrrrrr}
 \frac{7691}{684} & -\frac{9157}{1368} & \frac{11009}{4104} & -\frac{4561}{2052} & \frac{10501}{4104} & -\frac{2813}{2052} & -\frac{3673}{1368} & \frac{397}{1368} &
   \frac{286}{57} & \frac{10133}{2052} & -\frac{6077}{4104} & \frac{2209}{1026} & \frac{13}{216} & \frac{319}{4104} & -\frac{8237}{4104} & \frac{5005}{1368} &
   -\frac{443}{228} \\[2mm]
 -\frac{9157}{1368} & \frac{7691}{684} & -\frac{2813}{2052} & \frac{10501}{4104} & -\frac{4561}{2052} & \frac{11009}{4104} & \frac{286}{57} & \frac{397}{1368} &
   -\frac{3673}{1368} & -\frac{6077}{4104} & \frac{10133}{2052} & -\frac{8237}{4104} & \frac{319}{4104} & \frac{13}{216} & \frac{2209}{1026} & -\frac{443}{228} &
   \frac{5005}{1368} \\[2mm]
 -\frac{2345}{1368} & -\frac{2099}{171} & \frac{1132}{513} & \frac{9409}{4104} & \frac{2549}{513} & -\frac{12763}{4104} & -\frac{1939}{342} & \frac{1603}{1368} &
   -\frac{1033}{456} & -\frac{5915}{4104} & -\frac{5657}{1026} & \frac{20095}{4104} & \frac{13297}{4104} & \frac{11569}{4104} & \frac{1601}{2052} & \frac{47}{171} &
   -\frac{1429}{456} \\[2mm]
 \frac{1571}{342} & -\frac{16343}{1368} & \frac{11663}{4104} & \frac{1183}{1026} & \frac{20803}{4104} & -\frac{6641}{2052} & -\frac{2281}{456} & \frac{2275}{1368} &
   \frac{449}{342} & \frac{2411}{2052} & -\frac{20279}{4104} & \frac{101}{27} & \frac{7321}{4104} & \frac{5665}{4104} & -\frac{5789}{4104} & \frac{1099}{456} &
   -\frac{523}{171} \\[2mm]
 -\frac{16343}{1368} & \frac{1571}{342} & -\frac{6641}{2052} & \frac{20803}{4104} & \frac{1183}{1026} & \frac{11663}{4104} & \frac{449}{342} & \frac{2275}{1368} &
   -\frac{2281}{456} & -\frac{20279}{4104} & \frac{2411}{2052} & -\frac{5789}{4104} & \frac{5665}{4104} & \frac{7321}{4104} & \frac{101}{27} & -\frac{523}{171} &
   \frac{1099}{456} \\[2mm]
 -\frac{2099}{171} & -\frac{2345}{1368} & -\frac{12763}{4104} & \frac{2549}{513} & \frac{9409}{4104} & \frac{1132}{513} & -\frac{1033}{456} & \frac{1603}{1368} &
   -\frac{1939}{342} & -\frac{5657}{1026} & -\frac{5915}{4104} & \frac{1601}{2052} & \frac{11569}{4104} & \frac{13297}{4104} & \frac{20095}{4104} &
   -\frac{1429}{456} & \frac{47}{171} \\[2mm]
 \frac{4319}{684} & -\frac{2377}{684} & -\frac{20821}{5130} & -\frac{8057}{2052} & \frac{11303}{10260} & \frac{17159}{5130} & -\frac{1097}{342} & -\frac{2771}{3420}
   & \frac{14287}{3420} & \frac{21263}{10260} & -\frac{953}{2052} & -\frac{6568}{2565} & -\frac{4877}{2052} & \frac{7367}{10260} & \frac{1231}{2052} &
   \frac{3871}{3420} & -\frac{3107}{3420} \\[2mm]
 \frac{8473}{3420} & \frac{8473}{3420} & \frac{1211}{5130} & -\frac{16903}{10260} & -\frac{16903}{10260} & \frac{1211}{5130} & \frac{4241}{1710} & \frac{673}{855} &
   \frac{4241}{1710} & -\frac{1849}{10260} & -\frac{1849}{10260} & -\frac{5416}{2565} & -\frac{24079}{10260} & -\frac{24079}{10260} & -\frac{5416}{2565} &
   \frac{53}{171} & \frac{53}{171} \\[2mm]
 -\frac{2377}{684} & \frac{4319}{684} & \frac{17159}{5130} & \frac{11303}{10260} & -\frac{8057}{2052} & -\frac{20821}{5130} & \frac{14287}{3420} &
   -\frac{2771}{3420} & -\frac{1097}{342} & -\frac{953}{2052} & \frac{21263}{10260} & \frac{1231}{2052} & \frac{7367}{10260} & -\frac{4877}{2052} &
   -\frac{6568}{2565} & -\frac{3107}{3420} & \frac{3871}{3420} \\[2mm]
 -\frac{824}{285} & -\frac{1198}{171} & -\frac{13063}{5130} & \frac{2573}{2565} & \frac{7679}{2565} & \frac{5057}{5130} & -\frac{6337}{1710} & \frac{1663}{1140} &
   -\frac{1333}{3420} & -\frac{16613}{5130} & -\frac{8501}{2052} & -\frac{61}{135} & -\frac{1081}{10260} & \frac{2455}{2052} & \frac{2041}{2052} & -\frac{761}{1140}
   & -\frac{5467}{3420} \\[2mm]
 -\frac{1198}{171} & -\frac{824}{285} & \frac{5057}{5130} & \frac{7679}{2565} & \frac{2573}{2565} & -\frac{13063}{5130} & -\frac{1333}{3420} & \frac{1663}{1140} &
   -\frac{6337}{1710} & -\frac{8501}{2052} & -\frac{16613}{5130} & \frac{2041}{2052} & \frac{2455}{2052} & -\frac{1081}{10260} & -\frac{61}{135} &
   -\frac{5467}{3420} & -\frac{761}{1140} \\[2mm]
 -\frac{22037}{2280} & \frac{119207}{6840} & -\frac{3841}{20520} & \frac{66763}{20520} & -\frac{122123}{20520} & \frac{26969}{20520} & \frac{63959}{6840} &
   \frac{8}{19} & -\frac{4843}{1368} & -\frac{68621}{20520} & \frac{24209}{4104} & -\frac{78761}{20520} & -\frac{2923}{2565} & -\frac{8251}{2565} &
   -\frac{139}{4104} & -\frac{2169}{760} & \frac{6557}{1368} \\[2mm]
 \frac{1829}{1710} & \frac{1931}{380} & -\frac{2605}{2052} & -\frac{8411}{5130} & -\frac{30991}{10260} & \frac{734}{513} & \frac{674}{855} & -\frac{2713}{1140} &
   -\frac{139}{684} & \frac{6062}{2565} & \frac{40547}{10260} & -\frac{1303}{2565} & \frac{35}{54} & \frac{1259}{10260} & \frac{3188}{2565} & -\frac{115}{342} &
   \frac{1193}{1140} \\[2mm]
 \frac{1931}{380} & \frac{1829}{1710} & \frac{734}{513} & -\frac{30991}{10260} & -\frac{8411}{5130} & -\frac{2605}{2052} & -\frac{139}{684} & -\frac{2713}{1140} &
   \frac{674}{855} & \frac{40547}{10260} & \frac{6062}{2565} & \frac{3188}{2565} & \frac{1259}{10260} & \frac{35}{54} & -\frac{1303}{2565} & \frac{1193}{1140} &
   -\frac{115}{342} \\[2mm]
 \frac{119207}{6840} & -\frac{22037}{2280} & \frac{26969}{20520} & -\frac{122123}{20520} & \frac{66763}{20520} & -\frac{3841}{20520} & -\frac{4843}{1368} &
   \frac{8}{19} & \frac{63959}{6840} & \frac{24209}{4104} & -\frac{68621}{20520} & -\frac{139}{4104} & -\frac{8251}{2565} & -\frac{2923}{2565} &
   -\frac{78761}{20520} & \frac{6557}{1368} & -\frac{2169}{760} \\[2mm]
 -\frac{2081}{570} & \frac{125431}{6840} & -\frac{31475}{4104} & -\frac{3703}{1026} & -\frac{121739}{20520} & \frac{19553}{2052} & \frac{13529}{2280} &
   -\frac{23827}{6840} & \frac{619}{1710} & \frac{4469}{10260} & \frac{193087}{20520} & -\frac{17584}{2565} & -\frac{79721}{20520} & \frac{9343}{20520} &
   \frac{106477}{20520} & -\frac{21281}{6840} & \frac{1487}{285} \\[2mm]
 \frac{125431}{6840} & -\frac{2081}{570} & \frac{19553}{2052} & -\frac{121739}{20520} & -\frac{3703}{1026} & -\frac{31475}{4104} & \frac{619}{1710} &
   -\frac{23827}{6840} & \frac{13529}{2280} & \frac{193087}{20520} & \frac{4469}{10260} & \frac{106477}{20520} & \frac{9343}{20520} & -\frac{79721}{20520} &
   -\frac{17584}{2565} & \frac{1487}{285} & -\frac{21281}{6840}
\end{array}
\right)
\end{eqnarray*}
\end{landscape}

\begin{landscape}
\tiny
\begin{eqnarray*}
A_3 & = &
\left(
\begin{array}{rrrrrrrrrrrrrrrrr}
 \frac{13997}{1026} & \frac{15049}{4104} & \frac{7531}{4104} & -\frac{3062}{513} & -\frac{4225}{1368} & \frac{89}{342} & \frac{785}{456} & -\frac{3851}{1368} &
   \frac{120}{19} & \frac{7957}{1026} & \frac{15553}{4104} & -\frac{17}{54} & -\frac{449}{216} & -\frac{8171}{4104} & -\frac{9767}{4104} & \frac{3911}{1368} &
   \frac{47}{171} \\[2mm]
 \frac{22307}{4104} & \frac{54187}{4104} & \frac{8251}{4104} & -\frac{15943}{4104} & -\frac{10201}{1368} & \frac{131}{456} & \frac{8869}{1368} & -\frac{2399}{684} &
   \frac{3299}{1368} & \frac{18013}{4104} & \frac{27103}{4104} & -\frac{6887}{4104} & -\frac{4225}{2052} & -\frac{1870}{513} & -\frac{7727}{4104} & \frac{887}{1368}
   & \frac{1099}{456} \\[2mm]
 -\frac{1319}{2052} & -\frac{41957}{4104} & \frac{14623}{4104} & \frac{4501}{2052} & \frac{4577}{1368} & -\frac{863}{171} & -\frac{677}{152} & \frac{1525}{1368} &
   -\frac{42}{19} & -\frac{391}{513} & -\frac{19439}{4104} & \frac{4303}{1026} & \frac{12691}{4104} & \frac{4825}{4104} & -\frac{5069}{4104} & \frac{887}{1368} &
   -\frac{523}{171} \\[2mm]
 \frac{34169}{4104} & -\frac{43913}{4104} & \frac{10537}{4104} & -\frac{2677}{4104} & \frac{4643}{1368} & -\frac{1283}{456} & -\frac{6971}{1368} & \frac{409}{684} &
   \frac{209}{72} & \frac{14935}{4104} & -\frac{14153}{4104} & \frac{14425}{4104} & \frac{4181}{2052} & \frac{1547}{2052} & -\frac{7043}{4104} & \frac{3911}{1368} &
   -\frac{1429}{456} \\[2mm]
 -\frac{17155}{4104} & \frac{65953}{4104} & \frac{133}{216} & \frac{1415}{4104} & -\frac{9643}{1368} & \frac{1475}{1368} & \frac{3469}{456} & -\frac{1573}{684} &
   -\frac{339}{152} & \frac{1795}{4104} & \frac{30277}{4104} & -\frac{7655}{4104} & -\frac{299}{1026} & -\frac{2441}{1026} & \frac{1897}{4104} & -\frac{2347}{1368}
   & \frac{5005}{1368} \\[2mm]
 -\frac{16895}{2052} & -\frac{29369}{4104} & \frac{5893}{4104} & \frac{475}{108} & \frac{4133}{1368} & -\frac{773}{228} & -\frac{4919}{1368} & \frac{1643}{1368} &
   -\frac{980}{171} & -\frac{7639}{2052} & -\frac{15629}{4104} & \frac{7631}{2052} & \frac{15877}{4104} & \frac{9871}{4104} & \frac{5251}{4104} & -\frac{2347}{1368}
   & -\frac{443}{228} \\[2mm]
 -\frac{11509}{5130} & -\frac{1003}{10260} & -\frac{13723}{10260} & \frac{4391}{5130} & \frac{167}{180} & \frac{1411}{1140} & \frac{401}{285} & \frac{3181}{1710} &
   \frac{1061}{1710} & -\frac{7006}{2565} & -\frac{17719}{10260} & -\frac{20899}{10260} & -\frac{3826}{2565} & -\frac{8209}{10260} & -\frac{5029}{10260} &
   -\frac{106}{171} & \frac{53}{171} \\[2mm]
 \frac{2515}{1026} & -\frac{929}{270} & -\frac{7168}{2565} & -\frac{3935}{1026} & \frac{365}{342} & \frac{679}{380} & -\frac{599}{228} & -\frac{1189}{855} &
   \frac{5017}{3420} & \frac{3863}{10260} & -\frac{2255}{2052} & -\frac{3193}{2565} & -\frac{1199}{513} & \frac{13877}{10260} & \frac{965}{1026} & -\frac{191}{855}
   & -\frac{3107}{3420} \\[2mm]
 -\frac{5081}{1026} & \frac{889}{513} & \frac{991}{5130} & \frac{3661}{1026} & \frac{112}{171} & \frac{181}{1140} & \frac{1657}{1140} & \frac{1309}{684} &
   -\frac{1571}{855} & -\frac{7141}{2565} & -\frac{61}{540} & -\frac{977}{5130} & \frac{8447}{10260} & -\frac{1963}{10260} & \frac{1139}{2565} & -\frac{191}{855} &
   \frac{3871}{3420} \\[2mm]
 -\frac{163009}{20520} & -\frac{56513}{5130} & -\frac{2971}{2565} & \frac{66761}{20520} & \frac{5249}{1140} & -\frac{11933}{6840} & -\frac{2009}{380} &
   \frac{14827}{6840} & -\frac{1481}{360} & -\frac{118553}{20520} & -\frac{14051}{2052} & \frac{32347}{20520} & \frac{39097}{20520} & \frac{9701}{4104} &
   \frac{2269}{2052} & -\frac{1658}{855} & -\frac{2169}{760} \\[2mm]
 -\frac{42173}{10260} & -\frac{8518}{2565} & -\frac{533}{270} & \frac{21967}{10260} & \frac{784}{285} & \frac{1021}{1140} & -\frac{10439}{3420} & \frac{107}{228} &
   -\frac{466}{171} & -\frac{12499}{10260} & -\frac{25}{1026} & \frac{2999}{2565} & \frac{20039}{10260} & \frac{7787}{2565} & \frac{2981}{1026} & -\frac{2429}{3420}
   & -\frac{115}{342} \\[2mm]
 -\frac{3617}{2565} & \frac{11347}{5130} & -\frac{5068}{2565} & -\frac{151}{270} & -\frac{43}{285} & \frac{4129}{1710} & -\frac{109}{380} & -\frac{2231}{1710} &
   -\frac{3253}{3420} & \frac{1093}{2052} & \frac{22733}{10260} & \frac{160}{513} & \frac{415}{2052} & \frac{19073}{10260} & \frac{13849}{5130} & -\frac{2429}{3420}
   & \frac{1193}{1140} \\[2mm]
 -\frac{132199}{20520} & \frac{92413}{5130} & \frac{205}{513} & \frac{35951}{20520} & -\frac{7453}{1140} & \frac{3991}{2280} & \frac{6845}{684} & -\frac{179}{760} &
   -\frac{11701}{6840} & -\frac{40487}{20520} & \frac{65351}{10260} & -\frac{4313}{1080} & -\frac{2051}{1080} & -\frac{15133}{4104} & -\frac{5029}{10260} &
   -\frac{1658}{855} & \frac{6557}{1368} \\[2mm]
 \frac{33437}{5130} & -\frac{23953}{10260} & \frac{2863}{2052} & -\frac{3284}{2565} & \frac{509}{1140} & -\frac{523}{855} & \frac{47}{380} & \frac{3233}{3420} &
   \frac{688}{171} & \frac{9307}{5130} & -\frac{14779}{10260} & -\frac{2063}{5130} & -\frac{2609}{2052} & -\frac{4882}{2565} & -\frac{6526}{2565} & \frac{775}{342}
   & -\frac{761}{1140} \\[2mm]
 \frac{51557}{5130} & -\frac{51289}{10260} & \frac{3443}{10260} & -\frac{12344}{2565} & \frac{75}{76} & -\frac{21}{95} & -\frac{1189}{684} & -\frac{29}{57} &
   \frac{18733}{3420} & \frac{6691}{2052} & -\frac{20527}{10260} & -\frac{793}{1026} & -\frac{13943}{5130} & -\frac{9703}{10260} & -\frac{26377}{10260} &
   \frac{775}{342} & -\frac{5467}{3420} \\[2mm]
 -\frac{21365}{4104} & \frac{44899}{20520} & -\frac{7219}{4104} & \frac{403}{216} & -\frac{193}{1368} & \frac{9313}{6840} & \frac{13183}{6840} & \frac{6211}{3420} &
   -\frac{4751}{6840} & -\frac{84623}{20520} & -\frac{22907}{20520} & -\frac{54053}{20520} & -\frac{7481}{5130} & -\frac{12139}{10260} & -\frac{4067}{20520} &
   -\frac{7129}{6840} & \frac{7027}{6840} \\[2mm]
 -\frac{79909}{20520} & -\frac{30997}{4104} & -\frac{8299}{4104} & \frac{78761}{20520} & \frac{42259}{6840} & \frac{2113}{1368} & -\frac{29269}{6840} &
   \frac{1859}{570} & -\frac{2611}{2280} & -\frac{52403}{20520} & -\frac{77327}{20520} & \frac{19687}{20520} & \frac{8377}{5130} & \frac{32657}{10260} &
   \frac{48463}{20520} & -\frac{637}{6840} & -\frac{69}{760}
\end{array}
\right)
\\[3mm]
A_4 & = &
\left(
\begin{array}{rrrrrrrrrrrrrrrrr}
 \frac{1453}{152} & -\frac{11939}{1368} & \frac{307}{152} & -\frac{1985}{1368} & \frac{439}{152} & -\frac{2507}{1368} & -\frac{6071}{1368} & \frac{73}{342} &
   \frac{5017}{1368} & \frac{6221}{1368} & -\frac{2933}{1368} & \frac{3911}{1368} & \frac{989}{684} & \frac{517}{684} & -\frac{2003}{1368} & \frac{463}{152} &
   -\frac{385}{152} \\[2mm]
 -\frac{11939}{1368} & \frac{1453}{152} & -\frac{2507}{1368} & \frac{439}{152} & -\frac{1985}{1368} & \frac{307}{152} & \frac{5017}{1368} & \frac{73}{342} &
   -\frac{6071}{1368} & -\frac{2933}{1368} & \frac{6221}{1368} & -\frac{2003}{1368} & \frac{517}{684} & \frac{989}{684} & \frac{3911}{1368} & -\frac{385}{152} &
   \frac{463}{152} \\[2mm]
 -\frac{743}{76} & -\frac{5315}{1368} & -\frac{21}{8} & \frac{2693}{684} & \frac{385}{152} & \frac{1091}{684} & -\frac{4721}{1368} & \frac{859}{1368} &
   -\frac{3319}{684} & -\frac{737}{171} & -\frac{2705}{1368} & \frac{1027}{684} & \frac{3997}{1368} & \frac{4613}{1368} & \frac{6145}{1368} & -\frac{385}{152} &
   -\frac{39}{76} \\[2mm]
 \frac{9425}{684} & -\frac{5}{152} & \frac{1921}{1368} & -\frac{431}{76} & -\frac{95}{72} & \frac{27}{76} & -\frac{545}{1368} & -\frac{3317}{1368} & \frac{227}{36}
   & \frac{1289}{171} & \frac{3209}{1368} & \frac{7}{18} & -\frac{2113}{1368} & -\frac{1301}{1368} & -\frac{2431}{1368} & \frac{463}{152} & -\frac{39}{76} \\[2mm]
 -\frac{5}{152} & \frac{9425}{684} & \frac{27}{76} & -\frac{95}{72} & -\frac{431}{76} & \frac{1921}{1368} & \frac{227}{36} & -\frac{3317}{1368} & -\frac{545}{1368}
   & \frac{3209}{1368} & \frac{1289}{171} & -\frac{2431}{1368} & -\frac{1301}{1368} & -\frac{2113}{1368} & \frac{7}{18} & -\frac{39}{76} & \frac{463}{152} \\[2mm]
 -\frac{5315}{1368} & -\frac{743}{76} & \frac{1091}{684} & \frac{385}{152} & \frac{2693}{684} & -\frac{21}{8} & -\frac{3319}{684} & \frac{859}{1368} &
   -\frac{4721}{1368} & -\frac{2705}{1368} & -\frac{737}{171} & \frac{6145}{1368} & \frac{4613}{1368} & \frac{3997}{1368} & \frac{1027}{684} & -\frac{39}{76} &
   -\frac{385}{152} \\[2mm]
 \frac{59581}{3420} & -\frac{5923}{684} & \frac{2881}{3420} & -\frac{21659}{3420} & \frac{1937}{684} & \frac{1561}{3420} & -\frac{5729}{1710} & -\frac{143}{1710} &
   \frac{3221}{342} & \frac{469}{76} & -\frac{1031}{380} & -\frac{77}{228} & -\frac{3797}{1140} & -\frac{241}{228} & -\frac{3913}{1140} & \frac{1354}{285} &
   -\frac{677}{285} \\[2mm]
 -\frac{6187}{684} & -\frac{6187}{684} & -\frac{5411}{3420} & \frac{2201}{684} & \frac{2201}{684} & -\frac{5411}{3420} & -\frac{8267}{1710} & \frac{479}{342} &
   -\frac{8267}{1710} & -\frac{2207}{380} & -\frac{2207}{380} & \frac{1619}{1140} & \frac{2323}{1140} & \frac{2323}{1140} & \frac{1619}{1140} & -\frac{677}{285} &
   -\frac{677}{285} \\[2mm]
 -\frac{5923}{684} & \frac{59581}{3420} & \frac{1561}{3420} & \frac{1937}{684} & -\frac{21659}{3420} & \frac{2881}{3420} & \frac{3221}{342} & -\frac{143}{1710} &
   -\frac{5729}{1710} & -\frac{1031}{380} & \frac{469}{76} & -\frac{3913}{1140} & -\frac{241}{228} & -\frac{3797}{1140} & -\frac{77}{228} & -\frac{677}{285} &
   \frac{1354}{285} \\[2mm]
 \frac{37273}{6840} & -\frac{1471}{380} & -\frac{9607}{3420} & -\frac{11699}{2280} & \frac{56}{45} & \frac{6017}{2280} & -\frac{889}{342} & -\frac{2113}{1368} &
   \frac{1237}{360} & \frac{1931}{1368} & -\frac{4187}{3420} & -\frac{115}{72} & -\frac{21517}{6840} & \frac{7879}{6840} & \frac{1289}{1710} & \frac{26}{57} &
   -\frac{2569}{2280} \\[2mm]
 -\frac{1471}{380} & \frac{37273}{6840} & \frac{6017}{2280} & \frac{56}{45} & -\frac{11699}{2280} & -\frac{9607}{3420} & \frac{1237}{360} & -\frac{2113}{1368} &
   -\frac{889}{342} & -\frac{4187}{3420} & \frac{1931}{1368} & \frac{1289}{1710} & \frac{7879}{6840} & -\frac{21517}{6840} & -\frac{115}{72} & -\frac{2569}{2280} &
   \frac{26}{57} \\[2mm]
 -\frac{30023}{6840} & \frac{12133}{2280} & \frac{24793}{6840} & \frac{4019}{2280} & -\frac{23801}{6840} & -\frac{187}{40} & \frac{26437}{6840} & -\frac{1183}{3420}
   & -\frac{25979}{6840} & -\frac{1507}{1368} & \frac{9743}{6840} & \frac{1307}{1368} & \frac{176}{171} & -\frac{1939}{855} & -\frac{18397}{6840} &
   -\frac{2569}{2280} & \frac{1529}{2280} \\[2mm]
 -\frac{6499}{2280} & -\frac{433}{3420} & \frac{145}{228} & \frac{21673}{6840} & \frac{187}{95} & \frac{161}{1368} & \frac{851}{855} & \frac{17263}{6840} &
   -\frac{865}{1368} & -\frac{15277}{6840} & -\frac{3169}{3420} & -\frac{667}{6840} & \frac{683}{1368} & \frac{569}{6840} & \frac{43}{1710} & \frac{26}{57} &
   \frac{1529}{2280} \\[2mm]
 -\frac{433}{3420} & -\frac{6499}{2280} & \frac{161}{1368} & \frac{187}{95} & \frac{21673}{6840} & \frac{145}{228} & -\frac{865}{1368} & \frac{17263}{6840} &
   \frac{851}{855} & -\frac{3169}{3420} & -\frac{15277}{6840} & \frac{43}{1710} & \frac{569}{6840} & \frac{683}{1368} & -\frac{667}{6840} & \frac{1529}{2280} &
   \frac{26}{57} \\[2mm]
 \frac{12133}{2280} & -\frac{30023}{6840} & -\frac{187}{40} & -\frac{23801}{6840} & \frac{4019}{2280} & \frac{24793}{6840} & -\frac{25979}{6840} &
   -\frac{1183}{3420} & \frac{26437}{6840} & \frac{9743}{6840} & -\frac{1507}{1368} & -\frac{18397}{6840} & -\frac{1939}{855} & \frac{176}{171} & \frac{1307}{1368}
   & \frac{1529}{2280} & -\frac{2569}{2280} \\[2mm]
 -\frac{99}{38} & \frac{32563}{6840} & -\frac{511}{152} & -\frac{200}{171} & -\frac{353}{760} & \frac{59}{18} & \frac{16873}{6840} & \frac{4117}{6840} &
   \frac{1997}{3420} & -\frac{7913}{3420} & \frac{6241}{6840} & -\frac{6919}{1710} & -\frac{23003}{6840} & -\frac{2851}{6840} & \frac{4291}{6840} &
   -\frac{2803}{2280} & \frac{1091}{570} \\[2mm]
 \frac{32563}{6840} & -\frac{99}{38} & \frac{59}{18} & -\frac{353}{760} & -\frac{200}{171} & -\frac{511}{152} & \frac{1997}{3420} & \frac{4117}{6840} &
   \frac{16873}{6840} & \frac{6241}{6840} & -\frac{7913}{3420} & \frac{4291}{6840} & -\frac{2851}{6840} & -\frac{23003}{6840} & -\frac{6919}{1710} &
   \frac{1091}{570} & -\frac{2803}{2280}
\end{array}
\right),
\end{eqnarray*}
\end{landscape}

\normalsize

the four matrices $B_1, B_2, B_3, B_4$:
\tiny
\begin{eqnarray*}
&& \hskip -7mm B_1 \!=\!
\left(
\begin{array}{rrrrrrrrrrr}
 -\frac{245695}{30276} & \frac{1024085}{121104} & -\frac{146761}{121104} & \frac{68083}{60552} & -\frac{138935}{40368} & \frac{32845}{30276} & \frac{152425}{40368}
   & -\frac{107765}{121104} & -\frac{59195}{15138} & -\frac{51605}{15138} & -\frac{12635}{7569} \\[1.9mm]
 -\frac{27385}{7569} & \frac{935519}{121104} & \frac{47249}{121104} & \frac{40939}{60552} & -\frac{46215}{13456} & \frac{6745}{30276} & \frac{150395}{40368} &
   -\frac{122555}{121104} & -\frac{56185}{30276} & -\frac{25345}{30276} & -\frac{5675}{7569} \\[1.9mm]
 \frac{2425}{841} & -\frac{200725}{13456} & -\frac{69251}{121104} & -\frac{49109}{60552} & \frac{575405}{121104} & -\frac{22085}{30276} & -\frac{967175}{121104} &
   \frac{1945}{13456} & \frac{9595}{7569} & -\frac{2165}{7569} & \frac{22400}{7569} \\[1.9mm]
 -\frac{111245}{10092} & -\frac{36485}{40368} & -\frac{159509}{121104} & \frac{246115}{60552} & \frac{183455}{121104} & -\frac{950}{7569} & -\frac{118385}{121104} &
   \frac{51305}{40368} & -\frac{169535}{30276} & -\frac{162035}{30276} & \frac{19955}{30276} \\[1.9mm]
 -\frac{825575}{30276} & \frac{593215}{121104} & -\frac{53651}{13456} & \frac{659713}{60552} & \frac{252649}{121104} & \frac{9925}{7569} & \frac{144065}{121104} &
   \frac{430865}{121104} & -\frac{33280}{2523} & -\frac{30625}{2523} & -\frac{2305}{5046} \\[1.9mm]
 -\frac{111785}{30276} & -\frac{616655}{121104} & -\frac{106391}{40368} & \frac{43459}{60552} & \frac{245095}{121104} & \frac{47515}{30276} & -\frac{421105}{121104}
   & \frac{31715}{121104} & -\frac{4355}{3364} & -\frac{8615}{3364} & \frac{305}{3364} \\[1.9mm]
 -\frac{191761}{6728} & \frac{76347}{26912} & -\frac{718467}{134560} & \frac{841397}{67280} & \frac{120045}{26912} & \frac{3143}{1682} & \frac{13213}{26912} &
   \frac{153665}{26912} & -\frac{90001}{6728} & -\frac{89053}{6728} & -\frac{1307}{841} \\[1.9mm]
 \frac{218579}{6728} & \frac{47955}{26912} & \frac{2017423}{403680} & -\frac{2658161}{201840} & -\frac{424769}{80736} & -\frac{2486}{2523} & \frac{217223}{80736} &
   -\frac{124463}{26912} & \frac{321065}{20184} & \frac{307733}{20184} & -\frac{2953}{2523} \\[1.9mm]
 -\frac{738031}{20184} & \frac{373229}{80736} & -\frac{1481557}{403680} & \frac{1203897}{67280} & \frac{342619}{80736} & -\frac{163}{1682} & \frac{121883}{80736} &
   \frac{556927}{80736} & -\frac{361817}{20184} & -\frac{326087}{20184} & \frac{4297}{10092} \\[1.9mm]
 \frac{1401}{29} & -\frac{3977}{232} & \frac{21023}{3480} & -\frac{32659}{1740} & \frac{1055}{696} & -\frac{359}{174} & -\frac{4775}{696} & -\frac{1145}{232} &
   \frac{2039}{87} & \frac{7421}{348} & \frac{205}{87} \\[1.9mm]
 -\frac{4781}{174} & \frac{13555}{696} & -\frac{12787}{3480} & \frac{19261}{1740} & -\frac{893}{232} & \frac{223}{87} & \frac{2075}{232} & \frac{1985}{696} &
   -\frac{1111}{87} & -\frac{949}{87} & -\frac{1313}{348}
\end{array}
\right)
\\[2mm]
&& \hskip -7mm B_2 \!=\!
\left(
\begin{array}{rrrrrrrrrrr}
 \frac{1650899}{121104} & -\frac{940283}{121104} & \frac{33895}{121104} & -\frac{499045}{121104} & \frac{293725}{121104} & \frac{60065}{121104} &
   -\frac{465355}{121104} & -\frac{5}{696} & \frac{842065}{121104} & \frac{220075}{40368} & \frac{18395}{40368} \\[1.9mm]
 -\frac{186367}{121104} & \frac{896983}{121104} & -\frac{231455}{121104} & -\frac{35335}{121104} & -\frac{169985}{121104} & \frac{325415}{121104} &
   \frac{384635}{121104} & -\frac{5}{696} & -\frac{7925}{121104} & \frac{5185}{40368} & -\frac{100795}{40368} \\[1.9mm]
 \frac{383841}{13456} & -\frac{1294319}{40368} & \frac{462383}{40368} & -\frac{272725}{40368} & \frac{344065}{40368} & -\frac{154785}{13456} & -\frac{165345}{13456}
   & \frac{235}{696} & \frac{413965}{40368} & \frac{420985}{40368} & \frac{426805}{40368} \\[1.9mm]
 -\frac{55715}{40368} & -\frac{103615}{13456} & -\frac{113975}{40368} & \frac{17775}{13456} & \frac{54145}{13456} & \frac{57715}{40368} & -\frac{184385}{40368} &
   \frac{1075}{696} & -\frac{695}{13456} & -\frac{69905}{40368} & -\frac{485}{40368} \\[1.9mm]
 -\frac{1918919}{121104} & \frac{819239}{121104} & -\frac{908875}{121104} & \frac{547705}{121104} & \frac{99575}{121104} & \frac{740095}{121104} &
   \frac{149515}{121104} & \frac{1075}{696} & -\frac{708925}{121104} & -\frac{94545}{13456} & -\frac{171875}{40368} \\[1.9mm]
 \frac{2091859}{121104} & -\frac{2520247}{121104} & \frac{686915}{121104} & -\frac{465245}{121104} & \frac{679265}{121104} & -\frac{692831}{121104} &
   -\frac{1055135}{121104} & \frac{235}{696} & \frac{808925}{121104} & \frac{81845}{13456} & \frac{248455}{40368} \\[1.9mm]
 -\frac{4879959}{134560} & \frac{3060847}{134560} & -\frac{366567}{26912} & \frac{263909}{26912} & -\frac{130045}{26912} & \frac{314715}{26912} &
   \frac{185503}{26912} & \frac{287}{464} & -\frac{375917}{26912} & -\frac{409385}{26912} & -\frac{235637}{26912} \\[1.9mm]
 \frac{1103137}{134560} & \frac{1160719}{134560} & \frac{486971}{80736} & -\frac{253961}{80736} & -\frac{454103}{80736} & -\frac{258103}{80736} &
   \frac{461357}{80736} & -\frac{1085}{464} & \frac{219041}{80736} & \frac{436517}{80736} & \frac{110665}{80736} \\[1.9mm]
 -\frac{27188873}{403680} & \frac{21731537}{403680} & -\frac{1852889}{80736} & \frac{1366043}{80736} & -\frac{964451}{80736} & \frac{1697333}{80736} &
   \frac{1588793}{80736} & \frac{287}{464} & -\frac{2160035}{80736} & -\frac{711797}{26912} & -\frac{466593}{26912} \\[1.9mm]
 \frac{37341}{1160} & -\frac{37341}{1160} & \frac{4121}{232} & -\frac{1487}{232} & \frac{1487}{232} & -\frac{4121}{232} & -\frac{2411}{232} & 0 & \frac{2411}{232} &
   \frac{2851}{232} & \frac{3075}{232} \\[1.9mm]
 -\frac{178027}{3480} & \frac{178027}{3480} & -\frac{13987}{696} & \frac{8521}{696} & -\frac{8521}{696} & \frac{13987}{696} & \frac{13207}{696} & 0 &
   -\frac{13207}{696} & -\frac{4417}{232} & -\frac{3887}{232}
\end{array}
\right)
\\[2mm]
&& \hskip -7mm B_3 \!=\!
\left(
\begin{array}{rrrrrrrrrrr}
 \frac{661279}{121104} & \frac{679825}{30276} & \frac{815}{7569} & -\frac{192245}{40368} & -\frac{555283}{60552} & \frac{421333}{121104} & \frac{338935}{30276} &
   -\frac{441205}{121104} & \frac{147145}{40368} & \frac{650225}{121104} & -\frac{208325}{40368} \\[1.9mm]
 \frac{749845}{121104} & \frac{271835}{15138} & \frac{7340}{7569} & -\frac{192535}{40368} & -\frac{528139}{60552} & \frac{227323}{121104} & \frac{138365}{15138} &
   -\frac{426415}{121104} & \frac{49725}{13456} & \frac{558875}{121104} & -\frac{163375}{40368} \\[1.9mm]
 \frac{467615}{40368} & -\frac{84875}{5046} & \frac{29120}{7569} & -\frac{251765}{121104} & \frac{381797}{60552} & -\frac{553669}{121104} & -\frac{211075}{30276} &
   \frac{66955}{40368} & \frac{559145}{121104} & \frac{436655}{121104} & \frac{440815}{121104} \\[1.9mm]
 -\frac{148535}{40368} & \frac{29950}{2523} & \frac{12065}{30276} & \frac{203869}{121104} & -\frac{248029}{60552} & \frac{94085}{121104} & \frac{46190}{7569} &
   -\frac{16795}{40368} & -\frac{195685}{121104} & -\frac{15475}{121104} & -\frac{234455}{121104} \\[1.9mm]
 -\frac{1148275}{121104} & \frac{212810}{7569} & -\frac{31435}{30276} & \frac{134675}{121104} & -\frac{661627}{60552} & \frac{139145}{40368} & \frac{46065}{3364} &
   -\frac{327335}{121104} & -\frac{458135}{121104} & -\frac{165985}{121104} & -\frac{661045}{121104} \\[1.9mm]
 \frac{212975}{121104} & -\frac{310165}{30276} & \frac{11720}{7569} & \frac{78545}{121104} & \frac{289229}{60552} & -\frac{101249}{40368} & -\frac{11125}{2523} &
   \frac{186655}{121104} & \frac{13075}{121104} & -\frac{955}{121104} & \frac{308285}{121104} \\[1.9mm]
 -\frac{779423}{26912} & \frac{24497}{3364} & -\frac{30581}{6728} & \frac{300223}{26912} & \frac{8593}{67280} & \frac{303419}{134560} & \frac{15529}{6728} &
   \frac{62849}{26912} & -\frac{378439}{26912} & -\frac{333781}{26912} & -\frac{6733}{26912} \\[1.9mm]
 \frac{285081}{26912} & -\frac{81551}{3364} & -\frac{3631}{20184} & -\frac{287819}{80736} & \frac{1697971}{201840} & -\frac{734231}{403680} & -\frac{247579}{20184}
   & \frac{30097}{26912} & \frac{393923}{80736} & \frac{151169}{80736} & \frac{312241}{80736} \\[1.9mm]
 -\frac{2482457}{80736} & \frac{154865}{10092} & -\frac{17357}{6728} & \frac{918185}{80736} & -\frac{353907}{67280} & \frac{236413}{403680} & \frac{138401}{20184} &
   \frac{92615}{80736} & -\frac{1217561}{80736} & -\frac{359881}{26912} & -\frac{111259}{80736} \\[1.9mm]
 \frac{3977}{232} & -\frac{1401}{29} & \frac{359}{174} & -\frac{1055}{696} & \frac{32659}{1740} & -\frac{21023}{3480} & -\frac{2039}{87} & \frac{1145}{232} &
   \frac{4775}{696} & \frac{1919}{696} & \frac{6577}{696} \\[1.9mm]
 -\frac{13555}{696} & \frac{4781}{174} & -\frac{223}{87} & \frac{893}{232} & -\frac{19261}{1740} & \frac{12787}{3480} & \frac{1111}{87} & -\frac{1985}{696} &
   -\frac{2075}{232} & -\frac{4355}{696} & -\frac{1013}{232}
\end{array}
\right)
\\[2mm]
&& \hskip -7mm B_4 \!=\!
\left(
\begin{array}{rrrrrrrrrrr}
 \frac{79195}{15138} & -\frac{71365}{15138} & \frac{8855}{7569} & \frac{10381}{30276} & \frac{13805}{30276} & -\frac{38465}{30276} & -\frac{1775}{841} &
   \frac{35}{58} & \frac{12535}{5046} & \frac{47665}{30276} & \frac{7315}{7569} \\[1.9mm]
 -\frac{84595}{7569} & \frac{88510}{7569} & -\frac{92035}{30276} & \frac{116695}{30276} & -\frac{92509}{30276} & \frac{44495}{15138} & \frac{8195}{1682} &
   \frac{35}{58} & -\frac{11350}{2523} & -\frac{31910}{7569} & -\frac{99065}{30276} \\[1.9mm]
 -\frac{14795}{5046} & -\frac{57415}{5046} & -\frac{43595}{10092} & \frac{25265}{10092} & \frac{44915}{10092} & \frac{12851}{5046} & -\frac{24885}{3364} &
   \frac{95}{58} & -\frac{1615}{10092} & -\frac{5035}{1682} & \frac{305}{841} \\[1.9mm]
 \frac{10687}{1682} & \frac{24185}{3364} & \frac{7235}{10092} & -\frac{31055}{10092} & -\frac{8285}{1682} & \frac{1525}{1682} & \frac{12265}{3364} & -\frac{65}{29}
   & \frac{35705}{10092} & \frac{17335}{5046} & -\frac{17425}{10092} \\[1.9mm]
 -\frac{183155}{30276} & \frac{296593}{15138} & -\frac{9085}{15138} & \frac{9895}{15138} & -\frac{262085}{30276} & \frac{67325}{30276} & \frac{31985}{3364} &
   -\frac{65}{29} & -\frac{23455}{10092} & -\frac{32725}{30276} & -\frac{57965}{15138} \\[1.9mm]
 \frac{37685}{15138} & -\frac{254315}{15138} & -\frac{157}{15138} & \frac{39545}{30276} & \frac{170995}{30276} & -\frac{53365}{30276} & -\frac{28655}{3364} &
   \frac{95}{58} & \frac{9695}{10092} & -\frac{13415}{15138} & \frac{25220}{7569} \\[1.9mm]
 \frac{49977}{6728} & \frac{4762}{841} & \frac{10645}{1682} & -\frac{15721}{3364} & -\frac{32851}{6728} & -\frac{29733}{6728} & \frac{33881}{6728} &
   -\frac{773}{232} & \frac{1409}{841} & \frac{29543}{6728} & \frac{13711}{6728} \\[1.9mm]
 -\frac{6891}{6728} & -\frac{59673}{3364} & \frac{1263}{3364} & \frac{5295}{3364} & \frac{62751}{6728} & -\frac{17577}{6728} & -\frac{57239}{6728} & \frac{707}{232}
   & -\frac{2561}{1682} & -\frac{12589}{6728} & \frac{27207}{6728} \\[1.9mm]
 -\frac{515677}{20184} & \frac{97487}{2523} & -\frac{3967}{5046} & \frac{36125}{10092} & -\frac{265129}{20184} & \frac{54409}{20184} & \frac{130567}{6728} &
   -\frac{773}{232} & -\frac{42707}{3364} & -\frac{151259}{20184} & -\frac{115235}{20184} \\[1.9mm]
 \frac{2415}{116} & -\frac{2415}{116} & \frac{47}{116} & -\frac{869}{116} & \frac{869}{116} & -\frac{47}{116} & -\frac{1205}{116} & 0 & \frac{1205}{116} &
   \frac{455}{58} & \frac{291}{116} \\[1.9mm]
 -\frac{6433}{348} & \frac{6433}{348} & \frac{671}{348} & \frac{1795}{348} & -\frac{1795}{348} & -\frac{671}{348} & \frac{1205}{116} & 0 & -\frac{1205}{116} &
   -\frac{2237}{348} & -\frac{95}{87}
\end{array}
\right),
\end{eqnarray*}
\normalsize
and the four matrices $C_1, C_2, C_3, C_4$:
\tiny
\begin{eqnarray*}
C_1 & = &
\left(
\begin{array}{rrrrrr}
 1 & 0 & 0 & -\frac{1}{4} & 0 & 0 \\[1.9mm]
 0 & 0 & 0 & 0 & \frac{1}{4} & 0 \\[1.9mm]
 0 & 0 & 0 & 0 & 0 & \frac{1}{4} \\[1.9mm]
 0 & 0 & 0 & -\frac{1}{4} & 0 & 0 \\[1.9mm]
 0 & 1 & 0 & 0 & \frac{3}{4} & 0 \\[1.9mm]
 0 & 0 & 1 & 0 & 0 & \frac{3}{4}
\end{array}
\right)
\\[2mm]
C_2 & = &
\left(
\begin{array}{rrrrrr}
 0 & \frac{1}{4} & 0 & \frac{1}{12} & \frac{3}{16} & -\frac{1}{12} \\[1.9mm]
 1 & \frac{3}{4} & \frac{3}{4} & -\frac{23}{12} & -\frac{3}{16} & \frac{11}{48} \\[1.9mm]
 0 & 0 & -\frac{1}{4} & \frac{25}{12} & 0 & -\frac{25}{48} \\[1.9mm]
 0 & 0 & 0 & \frac{5}{12} & 0 & -\frac{1}{6} \\[1.9mm]
 0 & 0 & -1 & \frac{5}{3} & 1 & -\frac{5}{12} \\[1.9mm]
 0 & 0 & 0 & -\frac{7}{3} & 0 & \frac{1}{3}
\end{array}
\right)
\\[2mm]
C_3 & = &
\left(
\begin{array}{rrrrrr}
 0 & 0 & \frac{1}{4} & \frac{1}{12} & -\frac{1}{12} & \frac{3}{16} \\[1.9mm]
 0 & -\frac{1}{4} & 0 & \frac{25}{12} & -\frac{25}{48} & 0 \\[1.9mm]
 1 & \frac{3}{4} & \frac{3}{4} & -\frac{23}{12} & \frac{11}{48} & -\frac{3}{16} \\[1.9mm]
 0 & 0 & 0 & \frac{5}{12} & -\frac{1}{6} & 0 \\[1.9mm]
 0 & 0 & 0 & -\frac{7}{3} & \frac{1}{3} & 0 \\[1.9mm]
 0 & -1 & 0 & \frac{5}{3} & -\frac{5}{12} & 1
\end{array}
\right)
\\[2mm]
C_4 & = &
\left(
\begin{array}{rrrrrr}
 0 & -\frac{1}{12} & -\frac{1}{12} & -\frac{1}{4} & \frac{5}{48} & \frac{5}{48} \\[1.9mm]
 0 & \frac{5}{3} & -\frac{25}{12} & 0 & \frac{5}{3} & -\frac{103}{48} \\[1.9mm]
 0 & -\frac{25}{12} & \frac{5}{3} & 0 & -\frac{103}{48} & \frac{5}{3} \\[1.9mm]
 1 & \frac{1}{3} & \frac{1}{3} & -\frac{1}{4} & \frac{1}{3} & \frac{1}{3} \\[1.9mm]
 0 & -\frac{5}{3} & \frac{7}{3} & 0 & -\frac{23}{12} & \frac{25}{12} \\[1.9mm]
 0 & \frac{7}{3} & -\frac{5}{3} & 0 & \frac{25}{12} & -\frac{23}{12}
\end{array}
\right).
\end{eqnarray*}
\normalsize

\subsection*{Acknowledgments}

The research of the first author is supported by Italian
INdAM-GNCS (Gruppo Nazionale di Calcolo Scientifico). Part
of this work has been done during the stay of V.P. at the
University of L'Aquila and GSSI (L'Aquila).
Many thanks to Peter Oswald (Bremen) for providing us the
local subdivision matrices of the generalized Butterfly
subdivision scheme.

\bigskip

\end{document}